\documentclass{amsart}
\usepackage[dvips]{epsfig}
\usepackage{graphicx}
\usepackage{latexsym}
\usepackage{amsmath}
\usepackage{amsthm}
\usepackage{amssymb}
\usepackage{mathrsfs}

\usepackage{caption}
\usepackage{subcaption}

\usepackage{color}
\usepackage[utf8]{inputenc}

\newtheorem{thm}{Theorem}[section]
\newtheorem{theo}[thm]{Theorem}

\newtheorem{lem}[thm]{Lemma}
\newtheorem{cor}[thm]{Corollary}
\newtheorem{defi}[thm]{Definition}
\newtheorem{hyp}[thm]{Assumption}

\newtheorem{prop}[thm]{Proposition}
\theoremstyle{remark}

\newtheorem{rem}[thm]{Remark}

\newcommand{\vip}{\vskip.2cm}
\newcommand{\COMMENTAIRE}[1]{}

\newcommand{\field}[1]{\mathbb{#1}}
\newcommand{\EE}{\field{E}}

\newcommand{\GG}{\field{G}}

\newcommand{\II}{\field{I}}
\newcommand{\JJ}{\field{J}}

\newcommand{\NN}{\field{N}}
\newcommand{\PP}{\field{P}}

\newcommand{\RR}{\field{R}}
\newcommand{\TT}{\field{T}}

\newcommand{\Bb}{{\mathcal B}}
\newcommand{\Cc}{{\mathcal C}}

\newcommand{\Ff}{{\mathcal F}}

\newcommand{\Hh}{{\mathcal H}}
\newcommand{\Ll}{{\mathcal L}}

\newcommand{\Pp}{{\mathcal P}}

\newcommand{\Qq}{{\mathcal Q}}

\newcommand{\rd}{{\rm d}}

\newcommand{\bF}{{\mathfrak f}}

\newcommand{\cb}{{\mathcal B}}

\newcommand{\cf}{{\mathcal F}}

\newcommand{\cp}{{\mathcal P}}
\newcommand{\cq}{{\mathcal Q}}
\newcommand{\crr}{{\mathcal R}}

\newcommand{\cs}{{\mathscr S}}

\newcommand{\C}{{\mathbb C}}

\newcommand{\E}{{\mathbb E}}

\newcommand{\G}{\mathbb{G}}
\newcommand{\N}{{\mathbb N}}

\newcommand{\T}{\mathbb{T}}

\newcommand{\ind}{{\bf 1}}

\newcommand{\sot}{\otimes_{\rm sym}}
\newcommand{\ssub}{\Sigma^{\rm sub}}

\newcommand{\scrit}{\Sigma^{\rm crit}}

\newcommand{\inv}[1]{\mathop{\frac{1}{ #1}}\nolimits}

\newcommand{\reff}[1]{(\ref{#1})}

\begin{document}

\title[Moderate deviations principles for BMC]{Moderate deviation principles for bifurcating Markov chains: case of functions dependent of one variable}

\author{S. Val\`ere Bitseki Penda and Gorgui Gackou}
\address{S. Val\`ere Bitseki Penda, IMB, CNRS-UMR 5584, Universit\'e Bourgogne Franche-Comt\'e, 9 avenue Alain Savary, 21078 Dijon Cedex, France.}
\email{simeon-valere.bitseki-penda@u-bourgogne.fr}
\address{Gorgui Gackou, Laboratoire de Math\'ematiques Blaise Pascal, CNRS-UMR 6620, Universit\'e Clermont-Auvergne}
\email{gorgui.gackou@uca.fr}

\begin{abstract}
The main purpose of this article  is to establish moderate deviation principles for additive functionals of bifurcating Markov chains. Bifurcating Markov chains are a class of processes which are indexed by a regular binary tree. They can be seen as the models which represent the evolution of a trait along a population where each individual has two offsprings. Unlike the previous results of Bitseki, Djellout \& Guillin (2014), we consider here the case of functions which depend only on one variable. So, mainly inspired by the recent works of Bitseki \& Delmas (2020) about the  central limit theorem for general additive functionals of bifurcating Markov chains, we  give here  a moderate deviation principle for additive functionals of bifurcating Markov chains  when the functions depend on one variable. This work is done under  the uniform geometric ergodicity  and the uniform ergodic property based on the second spectral gap assumptions. The proofs of our results are based on martingale decomposition recently developed by Bitseki \& Delmas (2020) and on results of Dembo (1996), Djellout (2001) and Puhalski (1997).
\end{abstract}

\maketitle





\section{Introduction}

First, we give a general definition of a moderate deviation principles. Let $(Z_{n})_{n\geq 0}$ be a sequence of random variables with values in $S$ endowed with its Borel $\sigma$-field $\Bb(S)$ and let $(s_{n})_{n\geq 0}$ be a positive sequence that converges to $+\infty$. We assume that $Z_{n}/s_{n}$ converges in probability to 0 and that $Z_{n}/\sqrt{s_{n}}$ converges in distribution to a centered Gaussian law. Let $I: S \rightarrow\RR^{+}$ be a lower semicontinuous function, that is for all $c>0$ the sub-level set $\{x \in S, I(x)\leq c\}$ is a closed set. Such a function $I$ is called {\it rate function} and it is called {\it good rate function} if all its sub-level sets are compact sets. Let $(b_{n})_{n\geq 0}$ be a positive sequence such that $b_n \rightarrow + \infty $ and $b_{n}/\sqrt{s_{n}} \rightarrow 0$ as $n$ goes to $+\infty$.
\begin{defi}[Moderate deviation principle, MDP] \label{def:mdp}
	\quad\\
We say that $Z_{n}/(b_{n}\sqrt{s_{n}})$ satisfies a moderate deviation principle in $S$ with speed $b_{n}^{2}$ and the rate function~$I$ if, for any $A\in\Bb(S)$,
\begin{align*}
-\inf_{x\in A^{\circ}}I(x) \leq \liminf_{n\rightarrow\infty}\frac{1}{b_{n}^{2}} \log\PP\big(\frac{Z_{n}}{b_{n}\sqrt{s_{n}}}\in A\big)  \leq \limsup_{n\rightarrow\infty}\frac{1}{b_{n}^{2}} \log\PP\big(\frac{Z_{n}}{b_{n}\sqrt{s_{n}}}\in A\big) \leq -\inf_{x\in\bar{A}}I(x),
\end{align*}
where $A^{\circ}$ and $\bar{A}$ denote respectively the interior and the closure of $A$.
\end{defi}

\medskip

Bifurcating Markov chains (BMC, for short) are a class of stochastic processes indexed by regular binary tree.
They are appropriate for example  in the modeling of cell lineage data when each cell in one generation gives birth to two offspring in the next one.  Recently, they have received a great deal of attention because of the experiments of biologists on aging of Escherichia Coli (E. Coli , for  short). E. Coli is a rod-shaped bacterium which reproduces by dividing in two, thus producing two daughters: one of type 0 which has the old pole of the mother and the other of type 1 which has the new pole of the mother. The genealogy of the cells may be entirely described by a binary tree.
 To the best of our knowledge, the term bifurcating Markov chains appears for the first time in the works of Basawa and Zhou \cite{BZ04}. Thereafter, it was Guyon in \cite{Guyon} who had introduced and properly studied the theory of BMC. The first example of BMC, named bifurcating autoregressive process (BAR, for short), were introduced by Cowan and Staudte \cite{CS86} in order to study the mechanisms of cell division in Escherichia Coli.  Since this work of Cowan and Staudte, the BAR process has been widely studied in the literature and several extensions have been made. In particular, Guyon, in \cite{Guyon}, have used an extension of BAR process to get statistical evidence of aging in E.Coli. 

\medskip

In this paper, we are interested in moderate deviation principles( MDP, for short ) for additive functionals of bifurcating Markov chains.
The  MDP can be seen as an intermediate behavior between the central limit theorem   and large deviation. Usually, the MDP exhibits a simpler rate function inherited from the approximated Gaussian process, and holds for a larger class of dependent random variables than the large deviation principle.
 Unlike the results given in \cite{BDG14}, we treat here the case of functions which depends on one variable only. For this type of additive functionals, the martingale decomposition done in \cite{BDG14} is no longer valid. Indeed, as explained for e.g. in \cite{DM10} Remark 1.7, the error term on the last generation is not negligible. Note that recently, Bitseki and Delmas \cite{BD2020} have studied central limit theorem for additive functionals of bifurcating Markov chain. They have studied the case where the functions depend only on the trait of a single individual for BMC.  Bitseki and Delmas \cite{BD2020}	 observes three regimes (sub-critical, critical, super-critical), which correspond  to a competition between the reproducing rate (a mother has two daughters) and the ergodicity rate for the evolution of the trait along a lineage taken uniformly at random. This phenomenon already  appears in the works of Athreya \cite{athreya1969limit}. Here we investigate the moderate deviation principles for MBC depending only on one variable for the two cases: sub-critical  and  critical regimes. The super-critical regime, which require  another way of centering will be done in a future  work.
 
 \medskip

The rest of the paper is organized as follows. In Section \ref{sec:model}, we present the model of bifurcating Markov chains. In Section \ref{sec:nota}, we give some notations and the main assumptions for our results. In Section \ref{sec:main}, we set our main results: the sub-critical case in Section \ref{sec:main-S}  and the critical case in Section \ref{sec:main-Sc}. Section \ref{sec:proof-T1} is dedicated to the proof of the main result in sub-critical case and Section \ref{sec:proof-T2} is dedicated to the proof of the main result in Critical case. In Section \ref{sec:numerical}, we illustrate numerically our results.  Finally, in Section \ref{sec:appendix}, we give some useful results.  

\section{The model of bifurcating Markov chain}\label{sec:model}

\subsection{The regular binary tree associated to BMC models}
We denote by $\NN$ (resp. $\NN^{*}$) the space of (resp. positive) natural integers. We   set   $\T_0=\G_0=\{\emptyset\}$, $\G_k=\{0,1\}^k$  and $\T_k  =  \bigcup _{0  \leq r  \leq  k} \G_r$  for $k\in  \N^*$, and  $\T  =  \bigcup _{r\in  \N}  \G_r$. The set  $\G_k$ corresponds to the  $k$-th generation, $\T_k$ to the tree  up to the $k$-th generation, and $\T$ the complete binary  tree. One can see that the genealogy of the cells is entirely described by $\TT$ (each vertex of the tree designates an individual). For $i\in \T$, we denote by $|i|$ the generation of $i$ ($|i|=k$  if and only if $i\in \G_k$) and $iA=\{ij; j\in A\}$  for $A\subset \T$, where $ij$  is the concatenation of   the  two   sequences  $i,j\in   \T$,  with   the  convention   that $\emptyset i=i\emptyset=i$.

\subsection{The probability kernels associated to BMC models} 
\quad\\
 Let $(S,\cs)$ be a measurable space. For any $q \in \NN^{*}$, we denote by $\Bb(S^{q})$ (resp. $\Bb_{b}(S^{q})$, resp. $\Cc_{b}(S^{q})$) the space of  (resp. bounded, resp. bounded continuous ) $\RR\text{-}$valued measurable functions defined on $S^{q}$. For all $q \in \NN^{*}$, we set $\cs^{\otimes q} = \cs \otimes \ldots \otimes \cs$.  Let $\Pp$ be a probability kernel on $(S,\cs^{\otimes2})$, that is: $\Pp(\cdot  , A)$  is measurable  for all  $A\in \cs^{\otimes 2}$,  and $\Pp(x, \cdot)$ is  a probability measure on $(S^2,\cs^{\otimes 2})$ for all $x \in S$. For any $g\in \cb_b(S^3)$ and $h\in \cb_b(S^2)$,   we set for $x\in S$:
\begin{equation}\label{eq:Pg}
(\Pp g)(x)=\int_{S^2} g(x,y,z)\; \Pp(x,\rd y,\rd z) \quad \text{and} \quad (\Pp h)(x)=\int_{S^2} h(y,z)\; \Pp(x,\rd y,\rd z).
\end{equation}
We define $(\Pp g)$ (resp. $(\Pp h)$), or simply $Pg$ for $g\in \cb(S^3)$(resp. $\Pp h$ for $h\in \cb(S^2)$), as soon as the corresponding integral \reff{eq:Pg} is well defined, and we have  that $\Pp g$ and $\Pp h$ belong to $\cb(S)$. we denote by $\Pp_{0}$, $\Pp_{1}$ and $\Qq$ respectively the first and the second marginal of $\Pp$, and the mean of $\Pp_{0}$ and $\Pp_{1}$, that is, for all $x \in S$ and $B \in \mathcal{S}$  
\begin{equation*}\label{eq:P0P1Q}
\Pp_{0}(x,B) = \Pp(x,B\times S), \quad \Pp_{1}(x,B) = \Pp(x,S\times B) \quad \text{ and} \quad \Qq = \frac{(\Pp_{0} + \Pp_{1})}{2}.
\end{equation*}

\medskip 

Now let us give  a precise definition of bifurcating Markov chain. 
\begin{defi}[Bifurcating Markov Chains, see \cite{Guyon, BD2020}]
\quad\\
We say  a stochastic process indexed  by $\T$, $X=(X_i,  i\in \T)$, is a bifurcating Markov chain (BMC) on a measurable space $(S, \cs)$ with initial probability distribution  $\nu$ on $(S, \cs)$ and probability kernel $\cp$ on $S\times \cs^{\otimes 2}$ if:
\begin{itemize}
\item[-] (Initial  distribution.) The  random variable  $X_\emptyset$ is distributed as $\nu$.
\item[-] (Branching Markov property.) For  any  sequence   $(g_i, i\in \T)$ of functions belonging to $\cb_b(S^3)$ and  for all $k\geq 0$, we have
\begin{equation*}
\E\Big[\prod_{i\in \G_k} g_i(X_i,X_{i0},X_{i1}) |\sigma(X_j; j\in \T_k)\Big] 
=\prod_{i\in \G_k} \cp g_i(X_{i}).
\end{equation*}
\end{itemize}
\end{defi}
Following \cite{Guyon}, we introduce an  auxiliary Markov  chain $Y=(Y_n, n\in  \N) $  on $(S,\cs)$ with $Y_{0}=X_{1}$ and transition probability $\Qq$. The chain $(Y_{n}, n\in \mathbb{N})$ corresponds to a random lineage taken in the population. We  shall   write  $\E_x$   when  $X_\emptyset=x$ (\textit{i.e.}  the initial  distribution  $\nu$ is  the  Dirac mass  at $x\in S$).

\section{Notations and assumptions}\label{sec:nota}
For $f \in \Bb_{b}(S)$, we set $\|f\|_{\infty} = \sup\{|f(x)|, x \in S\}$. We will work with the following ergodic property.

\begin{hyp}\label{hyp:F2}
There exists a probability measure $\mu$ on $(S, \cs)$, a positive real number $M$  and $\alpha\in (0, 1)$ such that for all $f\in  \Bb_{b}(S)$: 
\begin{equation}\label{eq:geom-erg}
|\cq^{n}f - \langle \mu, f \rangle| \leq M \,  \alpha^{n} \|f\|_{\infty} \quad \text{for all  $n\in \N$.}
\end{equation}
\end{hyp}

We consider the stronger ergodic property based on a second spectral gap. 
\begin{hyp}\label{hyp:F3}
There exists a probability measure $\mu$ on $(S, \cs)$, a positive real number $M$, $\alpha\in (0, 1)$, a finite  non-empty set $J$ of indices, distinct complex eigenvalues $\{\alpha_j, \, j\in J\}$ of  the  operator  $\cq$ with  $|\alpha_j|=\alpha$,  non-zero  complex projectors $\{\crr_j, \, j\in J\}$  defined on $\C \Bb_{b}(S)$, the $\C$-vector space    spanned by        $\Bb_{b}(S)$,  such that $\crr_j\circ \crr_{j'}=\crr_{j'}\circ  \crr_{j}=0$ for all  $j\neq j'$ (so that  $\sum_{j\in J} \crr_j$ is  also a projector defined  on $\C \Bb_{b}(S)$) and a positive  sequence $(\beta_n, n\in \N)$ converging  to $0$, such that  for  all  $f\in  \Bb_{b}(S)$, with $\theta_j=\alpha_j/\alpha$:
\begin{equation}\label{eq:hyp-crit}
\Big|\cq^{n}( f) - \langle \mu, f \rangle - \alpha^n \sum_{j\in J}\theta_j^n\,  \crr_j (f) \Big| \leq M\, \beta_n \alpha^{n} \|f\|_{\infty} \quad \text{for all  $n\in \N$.}
\end{equation}
\end{hyp}
Without loss of generality, we shall assume that the sequence $(\beta_n,
n\in \N)$ in Assumption  \ref{hyp:F3} is non-increasing and bounded from above by 1. This assumption will be used when $\alpha = 1/\sqrt{2}.$ For $f \in \Bb_{b}(S)$, $\tilde{f}$ and $\hat{f}$ will denote the functions defined by:
\begin{equation}\label{eq:tilde-hat-f}
\tilde{f} = f - \langle f,\mu \rangle \quad \text{and} \quad \hat{f} = \tilde{f} - \alpha^{n} \sum_{j \in J} \theta_{j}^{n} \crr_{j}(f).
\end{equation}

\medskip

Let $\bF=(f_\ell, \ell\in \N)$ be a sequence of elements of $\Bb_{b}(S)$. We will assume in the sequel that 
\begin{equation}\label{eq:Bsupfell}
\sup_{\ell \in \NN} \{\|f_{\ell}\|_{\infty}\} = c_{\infty} < +\infty,
\end{equation}
in such a way that \eqref{eq:geom-erg} and \eqref{eq:hyp-crit} are uniformly satisfied by the sequence $\bF$.  We set for $n\in \N$ and $i\in \T_n$:
\begin{equation}\label{eq:def-NiF}
N_{n,i}(\bF)=\sum_{\ell=0}^{n-|i|} N_{n,i}^\ell(f_\ell)  = |\G_n|^{-1/2 }\sum_{\ell=0}^{n-|i|} M_{i\G_{n-|i|-\ell}}(\tilde f_\ell).
\end{equation}
We deduce that $ \sum_{i\in \G_k} N_{n,i}(\bF)=|\G_n|^{-1/2 }\sum_{\ell=0}^{n-k} M_{\G_{n-\ell}}(\tilde f_\ell)$ which gives for $k=0$ that 
\begin{equation*}\label{eq:def-NOf}
\boxed{N_{n, \emptyset}(\bF)= |\G_n|^{-1/2 }\sum_{\ell=0}^{n}M_{\G_{n-\ell}}(\tilde f_{\ell}).}
\end{equation*}
To study the asymptotics of $N_{n, \emptyset}(\bF)$, it is convenient to write for $n\geq k\geq 1$:
\begin{equation}\label{eq:nof-D}
N_{n, \emptyset}(\bF)= |\G_n|^{-1/2} \sum_{r=0}^{k-1} M_{\G_r}(\tilde{f}_{n-r}) + \sum_{i\in \G_k} N_{n,i}(\bF).
\end{equation} 

Asymptotic normality for $N_{n,\emptyset}(\bF)$ have been studied in \cite{BD2020}. Our aim in this paper is to complete this result by studying moderate deviation principles for $N_{n,\emptyset}(\bF)$. More precisely, given a sequence $(b_{n}, n \in \NN)$ such that:
\begin{equation*}
\lim_{n \rightarrow \infty} b_{n} = \infty \quad \text{and} \quad \lim_{n \rightarrow \infty} \frac{b_{n}}{\sqrt{|\GG_{n}|}} = 0,
\end{equation*}
our aim is to prove that $b_{n}^{-1}N_{n, \emptyset}(\bF)$ satisfies a moderate deviation principle with speed $b_{n}^{2}$ and rate function $I$ defined by
\begin{equation}\label{eq:mdp-rate}
I(x) = \sup_{\lambda \in \RR}\{\lambda x - \tfrac{1}{2} \lambda^{2} \Sigma(\bF)^{-1}\} = \begin{cases} \tfrac{1}{2} \Sigma(\bF)^{-1}x^{2} & \text{if $\Sigma(\bF) \neq 0$} \\ +\infty & \text{if $\Sigma(\bF) = 0$},  \end{cases}
\end{equation}
where
\begin{equation*}
\Sigma(\bF) = \begin{cases} \ssub(\bF)=\ssub_1(\bF)+ 2 \ssub_2(\bF) & \text{if $2\alpha^{2} < 1$} \\ \scrit (\bF)=\scrit _1(\bF)+ 2\scrit _2(\bF) & \text{if $2\alpha^{2} = 1,$} \end{cases}
\end{equation*}
with
\begin{align}
\label{eq:S1}
\ssub_1(\bF) &=\sum_{\ell\geq 0}  2^{-\ell} \, \langle \mu,   \tilde f_\ell^ 2\rangle + \sum_{\ell\geq 0, \, k\geq 0} 2^{k-\ell} \, \langle \mu, \cp\left((\cq^k \tilde f_\ell) \otimes^2\right)\rangle,\\
\label{eq:S2}
\ssub_2(\bF) &=\sum_{0\leq \ell< k} 2^{-\ell} \langle \mu,   \tilde f_k \cq^{k-\ell} \tilde f_\ell\rangle + \sum_{\substack{0\leq \ell< k\\ r\geq 0}} 2^{r-\ell} \langle \mu, \cp\left( \cq^r \tilde f_k \sot \cq^{k-\ell+r} \tilde f_\ell  \right)\rangle,\\
\label{eq:S1-crit}
\scrit_1(\bF)  &= \sum_{k\geq 0} 2^{-k} \langle \mu, \cp f_{k, k}^* \rangle = \sum_{k\geq 0} 2^{-k} \sum_{j\in J} \langle \mu, \cp(\crr_{j}(f_k) \sot \overline{\crr}_{j}(f_k)) \rangle,\\
\label{eq:S2-crit}
\scrit_2(\bF) &=   \sum_{0\leq \ell <k} 2^{-(k+\ell)/2} \langle \mu, \cp f_{k, \ell}^*  \rangle, 
\end{align}
and where for $k, \ell\in \N$:
\begin{equation}\label{eq:def-fkl*}
f_{k, \ell}^*= \sum_{j\in J}\,\, \theta_j^{\ell -k}\,  \crr_j
(f_k)\sot   \overline \crr_j (f_\ell).  
\end{equation}
More precisely, our aim is to prove that
\begin{align*}
-\inf_{x\in A^{\circ}}I(x) \leq \liminf_{n\rightarrow\infty}\frac{1}{b_{n}^{2}} \log\PP\big(b_{n}^{-1}N_{n, \emptyset}(\bF)\in A\big)  \leq \limsup_{n\rightarrow\infty}\frac{1}{b_{n}^{2}} \log\PP\big(b_{n}^{-1}N_{n, \emptyset}(\bF)\in A\big) \leq -\inf_{x\in\bar{A}}I(x),
\end{align*}
where $A^{\circ}$ and $\bar{A}$ denote respectively the interior and the closure of $A$. In particular, the latter asymptotic result implies that
\begin{equation*}
\lim_{n \rightarrow \infty} \frac{1}{b_{n}^{2}} \log \PP\left(\left|b_{n}^{-1} N_{n,\emptyset}(\bF)\right| > \delta \right) = -I(\delta) \quad \forall \delta > 0.
\end{equation*}
We note that $2\alpha^{2} < 1$ corresponds to the sub-critical regime and $2\alpha^{2} = 1$ to the critical regime. The super-critical regime, that is the case where  $2\alpha^{2} > 1$, is not treated in this paper. Indeed, for this case, another way to centered the functions is necessary to get moderate deviation principles. This will be done in a future work.

\begin{rem}\label{rem:NbF=Nf}
Let $f \in \Bb_{b}(S).$ If the sequence $\bF = (f_{\ell}, \ell \in \NN)$ is defined by: $f_{0} = f$ and $f_{\ell} = 0$ for all $\ell \geq 1$, then we have $N_{n,\emptyset}(\bF) = |\GG_{n}|^{-1/2} M_{\GG_{n}}(\tilde{f})$ and $\Sigma(\bF) = \Sigma_{\GG}(f)$, where
\begin{equation*}
\Sigma_{\GG}(f) = \begin{cases} \ssub_{\GG}(f) = \langle \mu,   \tilde f^ 2\rangle + \sum_{k\geq 0} 2^{k} \, \langle \mu, \cp\left(\cq^k \tilde f \otimes ^2\right)\rangle & \text{if $2\alpha^{2} < 1$} \\ \scrit_{\G}(f) = \sum_{j\in J} \langle \mu, \cp(\crr_{j}(f) \sot \overline{\crr}_{j}(f)) \rangle & \text{if $2\alpha^{2} = 1,$} \end{cases}
\end{equation*}
If the sequence $\bF = (f_{\ell}, \ell \in \NN)$ is defined by: $f_{\ell} = f$ for all $\ell \in \NN$, then we have $N_{n,\emptyset}(\bF) = |\GG_{n}|^{-1/2} M_{\TT_{n}}(\tilde{f}) = \sqrt{2 - 2^{-n}} \, |\TT_{n}|^{-1/2} \,  M_{\TT_{n}}(\tilde{f})$ and $\Sigma(\bF) = \Sigma_{\TT}(f)$, where
\begin{equation*}
\Sigma_{\TT}(f) = \begin{cases} \ssub_{\TT}(f) = \ssub_{\G}(f) +2  \ssub_{\TT,2}(f) & \text{if $2\alpha^{2} < 1$} \\ \scrit_{\TT}(f) = \scrit_{\G}(f) + 2\scrit_{\TT,2}(f) & \text{if $2\alpha^{2} = 1,$} \end{cases}
\end{equation*}
with
\begin{align*}
&\ssub_{\TT,2}(f) =    \sum_{k\geq 1} \langle \mu, \tilde  f \cq^{k} \tilde  f\rangle +  \sum_{\substack{k\geq 1 \\   r\geq  0}} 2^{r}   \langle   \mu, \cp\left(\cq^r  \tilde f \sot \cq^{r+k} \tilde f \right)\rangle, \\ 
&\scrit_{\TT,2}(f) =   \sum_{j\in J} \inv{\sqrt{2}\, \theta_j -1} \langle \mu, \cp(\crr_{j}(f) \sot \overline{\crr}_{j}(f)) \rangle .
\end{align*}
 
\end{rem}   

\medskip

\section{The main results}\label{sec:main}

\subsection{The sub-critical cases: $2\alpha^{2} < 1$}\label{sec:main-S}
\quad\\
In the sub-critical case, we consider a sequence $(b_{n}, n \in \NN)$ such that:
\begin{equation*}
\lim_{n \rightarrow \infty} b_{n} = \infty \quad \text{and} \quad \lim_{n \rightarrow \infty} \frac{b_{n}}{\sqrt{|\GG_{n}|}} = 0.
\end{equation*}
Then, we have the following result.
\begin{theo}\label{theo:mdp-SC}
Let  $X$  be  a  BMC   with  kernel  $\cp$  and  initial  distribution
$\nu$ such that Assumption \ref{hyp:F2} is in force  with $\alpha\in  (0, 1/\sqrt{2})$.  Let $\bF=(f_\ell, \ell\in  \N)$ be a sequence of elements of $\Bb_{b}(S)$ satisfying \eqref{eq:Bsupfell} and Assumption \ref{hyp:F2} uniformly. Then $b_{n}^{-1} N_{n,\emptyset}(\bF)$ satisfies a moderate deviation principle with speed $b_{n}^{2}$ and rate function $I$ defined in \eqref{eq:mdp-rate}.
\end{theo}

As a direct consequence of Remark \ref{rem:NbF=Nf} and Theorem \ref{theo:mdp-SC}, we have the following result.

\begin{cor}\label{cor:mdp-SC}
Let  $X$  be  a  BMC   with  kernel  $\cp$  and  initial  distribution
$\nu$ such that Assumption \ref{hyp:F2} is in force  with $\alpha\in  (0, 1/\sqrt{2})$. Let $f \in \Bb_{b}(S)$. Then $b_{n}^{-1} |\GG_{n}|^{-1/2} M_{\GG_{n}}(\tilde{f})$ and $b_{n}^{-1} |\TT_{n}|^{-1/2} M_{\TT_{n}}(\tilde{f})$ satisfy a moderate deviation principle with speed $b_{n}^{2}$ and rate function $I$ defined in \eqref{eq:mdp-rate}, with $\Sigma(\bF)$ replaced respectively by $\Sigma_{\GG}(f)$ and $\Sigma_{\TT}(f).$
\end{cor}

\subsection{The critical cases: $2\alpha^{2} = 1$}\label{sec:main-Sc}
\quad\\
In this  critical case, we consider a sequence $(b_{n}, n \in \NN)$ such that:
\begin{equation}\label{eq:vitesse-Sc}
\lim_{n \rightarrow \infty} b_{n} = \infty \quad \text{and} \quad \lim_{n \rightarrow \infty} \frac{b_{n}}{\sqrt{n|\GG_{n}|}} = 0.
\end{equation}
Then, we have the following result.
\begin{theo}\label{theo:mdp-C}
Let  $X$  be  a  BMC   with  kernel  $\cp$  and  initial  distribution
$\nu$ such that Assumption \ref{hyp:F3} is in force  with $\alpha = 1/\sqrt{2}$.  Let $\bF=(f_\ell, \ell\in  \N)$ be a sequence of elements of $\Bb_{b}(S)$ satisfying \eqref{eq:Bsupfell} and Assumption  \ref{hyp:F3} uniformly. Then $b_{n}^{-1} n^{-\frac{1}{2}} N_{n,\emptyset}(\bF)$ satisfies a moderate deviation principle with speed $b_{n}^{2}$ and rate function $I$ defined in \eqref{eq:mdp-rate}.
\end{theo}

As a direct consequence of Remark \ref{rem:NbF=Nf} and Theorem \ref{theo:mdp-C}, we have the following result.

\begin{cor}\label{cor:mdp-SC}
Let  $X$  be  a  BMC   with  kernel  $\cp$  and  initial  distribution
$\nu$ such that Assumption \ref{hyp:F3} is in force  with $\alpha =  1/\sqrt{2}$. Let $f \in \Bb_{b}(S)$. Then $b_{n}^{-1} (n|\GG_{n}|)^{-1/2} M_{\GG_{n}}(\tilde{f})$ and $b_{n}^{-1} (n|\TT_{n}|)^{-1/2} M_{\TT_{n}}(\tilde{f})$ satisfy a moderate deviation principle with speed $b_{n}^{2}$ and rate function $I$ defined in \eqref{eq:mdp-rate}, with $\Sigma(\bF)$ replaced respectively by $\Sigma_{\GG}(f)$ and $\Sigma_{\TT}(f).$
\end{cor}

\section{Proof of Theorem \ref{theo:mdp-SC}}\label{sec:proof-T1}

\subsection{A quick overview of our strategy}
\quad\\
Let $(p_n, n\in \N)$ be a non-decreasing sequence of elements of $\N^*$ such that:
\begin{equation}\label{eq:def-pn}
p_n< \frac{n}{2}.
\end{equation}
When there is no ambiguity, we write $p$ for $p_n$. 
\medskip

Let $i,j\in \T$. We write $i\preccurlyeq  j$ if $j\in i\T$. We denote by $i\wedge j$  the most recent  common ancestor of  $i$ and $j$,  which is defined  as   the  only   $u\in  \T$   such  that   if  $v\in   \T$  and $ v\preccurlyeq i$, $v \preccurlyeq j$  then $v \preccurlyeq u$. We also define the lexicographic order $i\leq j$ if either $i \preccurlyeq j$ or $v0  \preccurlyeq i$  and $v1  \preccurlyeq j$  for $v=i\wedge  j$.  Let $X=(X_i, i\in  \T)$ be  a $BMC$  with kernel  $\cp$ and  initial measure $\nu$. For $i\in \T$, we define the $\sigma$-field:
\begin{equation*}\label{eq:field-Fi}
\cf_{i}=\{X_u; u\in \T \text{ such that  $u\leq i$}\}.
\end{equation*}
By construction,  the $\sigma$-fields $(\cf_{i}; \, i\in \T)$ are nested
as $\cf_{i}\subset \cf_{j} $ for $i\leq  j$. 

\medskip

We define for $n\in \N$, $i\in \G_{n-p_n}$ and $\bF\in F^\N$ the martingale increments:
\begin{equation}\label{eq:def-DiF}
\Delta_{n,i}(\bF)= N_{n,i}(\bF) - \E\left[N_{n,i}(\bF) |\, \cf_i\right]\quad\text{and}\quad \Delta_n(\bF) = \sum_{i\in \G_{n-p_n}} \Delta_{n,i}(\bF).
\end{equation}
Thanks to \reff{eq:def-NiF}, we have:
\[
\sum_{i\in \G_{n-p_n}} N_{n, i}(\bF) = |\G_n|^{-1/2} \sum_{\ell=0}^{p_n}  M_{\G_{n-\ell}} (\tilde f_\ell) = |\G_n|^{-1/2} \sum_{k=n-p_n}^{n}  M_{\G_{k}} (\tilde f_{n-k}).
\]
Using the branching Markov property, and \eqref{eq:def-NiF}, we get for $i\in \G_{n-p_n}$:
\[ 
\E\left[N_{n,i}(\bF) |\, \cf_i\right] = \E\left[N_{n,i}(\bF) |\,
  X_i\right] = |\G_n|^{-1/2} \sum_{\ell=0}^{p_n} \E_{X_i}\left[M_{\G_{p_n-\ell}}(\tilde f_\ell)\right].
\] 
We deduce from \reff{eq:nof-D} with $k=n-p_n$ that:
\begin{equation}\label{eq:N=D+R}
N_{n, \emptyset}(\bF) = \Delta_n(\bF) + R_0(n)+R_1(n),
\end{equation}
with
\begin{equation*}\label{eq:reste01}
R_0(n)= |\G_n|^{-1/2}\, \sum_{k=0}^{n-p_n-1} M_{\G_k}(\tilde f_{n-k}) \quad\text{and}\quad R_1(n)= \sum_{i\in \G_{n-p_n}}\E\left[N_{n,i}(\bF) |\, \cf_i\right].
\end{equation*}

\medskip

Our goals will be achieved if we prove the following:

\begin{align}
&\label{eq:cvsupexpR0}
\forall \delta > 0, \quad \lim_{n \rightarrow \infty} \frac{1}{b_{n}^{2}} \log \PP(|b_{n}^{-1} R_{0}(n)| > \delta) = - \infty; \\
&\label{eq:cvsupexpR1}
\forall \delta > 0, \quad \lim_{n \rightarrow \infty} \frac{1}{b_{n}^{2}} \log \PP(|b_{n}^{-1} R_{1}(n)| > \delta) = - \infty; \\
&\label{eq:mdp-Dnf-SCr}
b_{n}^{-1}\Delta_{n}(\bF) \quad \text{satisfies a MDP on $S$ with speed $b_n^2$ and rate function $I$.}
\end{align}
Note that \eqref{eq:cvsupexpR0} and \eqref{eq:cvsupexpR1} mean that $R_{0}(n)$ and $R_{1}(n)$ are negligible in the sense of moderate deviations in such a way that from \eqref{eq:N=D+R}, $N_{n,\emptyset}(\bF)$ and $\Delta_{n}(\bF)$ satisfy the same moderate deviation principle (see Dembo and Zeitouni \cite{DZ1998}, chap. 4).


\subsection{Proof of \eqref{eq:cvsupexpR0}}\label{sec:supER0}
\quad\\
Using the Chernoff inequality, we have, for all $\lambda > 0$,
\begin{equation}\label{eq:cher-Scrit}
\PP(b_{n}^{-1} R_{0}(n) > \delta) \leq \exp(-\lambda b_{n} |\GG_{n}|^{1/2}) \, \EE\left[ \exp\left(\lambda \sum_{k = 0}^{n-p-1} M_{\GG_{k}}(\tilde{f}_{n-k})\right)\right].
\end{equation}
For all $\ell \in \{0, \ldots, n-p-1\}$, we set

\begin{equation*}
\II_{\ell} = \EE\left[ \exp\left(\lambda \sum_{k = 0}^{n - p - \ell - 2} M_{\GG_{k}}(\tilde{f}_{n-k})\right) \exp\left(\lambda M_{\GG_{n-p-\ell-1}}\left(\sum_{r=0}^{\ell} g_{p,r,\ell}\right)\right)\right],
\end{equation*}
where $g_{p,r,\ell} = 2^{r} \Qq^{r} \tilde{f}_{p+\ell+1-r}$, with the convention that an empty sum is zero. For all $\ell \in \{0, \ldots, n-p-2\}$, we have the following decomposition:
\begin{equation}\label{eq:DIIell}
\II_{\ell} = \EE\left[ \exp\left(\lambda \sum_{k = 0}^{n - p - \ell - 2} M_{\GG_{k}}(\tilde{f}_{n-k})\right)  \exp\left(\lambda M_{\GG_{n-p-\ell-2}} \left(\sum_{r=0}^{\ell} 2 \, \Qq(g_{p,r,\ell}) \right)\right) \JJ_{\ell}\right],
\end{equation}
where
\begin{equation*}
\JJ_{\ell} = \EE\left[\exp\left(\lambda \sum_{i \in \GG_{n-p-\ell-2}} \sum_{r = 0}^{\ell} \left(g_{p,r,\ell}(X_{i0}) + g_{p,r,\ell}(X_{i1}) - 2 \Qq (g_{p,r,\ell})(X_{i}) \right)\right) \Big| \Hh_{n-p-\ell-2}\right]. 
\end{equation*}
Using branching Markov property, we get
\begin{equation*}
\JJ_{\ell} = \prod_{i \in \GG_{n-p-\ell-2}} \EE_{X_{i}}\left[ \exp\left(\lambda \sum_{r = 0}^{\ell} \left(g_{p,r,\ell}(X_{i0}) + g_{p,r,\ell}(X_{i1}) - 2 \Qq (g_{p,r,\ell})(X_{i}) \right) \right)\right].
\end{equation*}
Using \eqref{eq:geom-erg} and \eqref{eq:Bsupfell}, we get
\begin{equation*}
\left|\sum_{r = 0}^{\ell} \left(g_{p,r,\ell}(X_{i0}) + g_{p,r,\ell}(X_{i1}) - 2 \Qq (g_{p,r,\ell})(X_{i})\right)\right| \leq 2M c_{\infty}\, \sum_{r=0}^{\ell} (2\, \alpha)^{r}.
\end{equation*}
Using Lemma \ref{lem:ABH} and the latter inequality, we get, for all $i \in \GG_{n - p - \ell - 2}$,
\begin{equation*}
\EE_{X_{i}}\left[ \exp\left(\lambda \sum_{r = 0}^{\ell} \left(g_{p,r,\ell}(X_{i0}) + g_{p,r,\ell}(X_{i1}) - 2 \Qq (g_{p,r,\ell})(X_{i}) \right) \right)\right] \leq \exp\left(2 \lambda^{2} M^{2} c_{\infty}^{2} (1 + \alpha)^{2} a_{\ell}^{2}\right),
\end{equation*}
with $a_{\ell} = \sum_{r=0}^{\ell} (2\alpha)^{r}$. The latter inequality implies that
\begin{equation}\label{eq:iJJell}
\JJ_{\ell} \leq \exp\left(2 \lambda^{2} M^{2} c_{\infty}^{2} (1 + \alpha)^{2} a_{\ell}^{2} |\GG_{n - p - \ell - 2}|\right).
\end{equation}
From \eqref{eq:DIIell} and \eqref{eq:iJJell}, it follows that
\begin{equation}\label{eq:recIIell}
\II_{\ell} \leq \exp\left(2 \lambda^{2} M^{2} c_{\infty}^{2} (1 + \alpha)^{2} a_{\ell}^{2} |\GG_{n - p - \ell - 2}|\right) \, \II_{\ell + 1}.
\end{equation}
Using the recurrence \eqref{eq:recIIell} for all $\ell \in \{0, \ldots, n-p-2\}$ for the first inequality, \eqref{eq:hyp-crit} and \eqref{eq:Bsupfell} for the second inequality, we are led to
\begin{multline*}
\EE\left[ \exp\left(\lambda \sum_{k = 0}^{n-p-1} M_{\GG_{k}}(\tilde{f}_{n-k})\right)\right] = \II_{0} \leq \exp\left(2 \lambda^{2} M^{2} c_{\infty}^{2} (1 + \alpha)^{2} \sum_{\ell = 0}^{n-p-2} a_{\ell}^{2} |\GG_{n-p-\ell-2}|\right) \, \II_{n-p-1} \\  \leq \exp\left(2 \lambda^{2} M^{2} c_{\infty}^{2} (1 + \alpha)^{2} \sum_{\ell = 0}^{n-p-2} a_{\ell}^{2} |\GG_{n-p-\ell-2}| + \lambda \, c_{\infty}  \, M \sum_{r = 0}^{n-p-1} (2\alpha)^{r+1}\right).
\end{multline*}

We have
\begin{equation*}
\sum_{\ell = 0}^{n-p-2} a_{\ell}^{2} |\GG_{n-p-\ell-2}| \leq \begin{cases} 6 \, |\GG_{n-p-1}| & \text{if $2\alpha \leq 1$} \\ \frac{2\alpha^{2}}{(2\alpha - 1)^{2}(1 - 2\alpha^{2})} \, |\GG_{n-p-1}| & \text{if $1 < 2 \alpha < \sqrt{2} $} \\ \frac{1}{(2\alpha - 1)^{2}}(n-p-1) |\GG_{n-p-1}| & \text{if $2\alpha^{2} = 1$} \\ \frac{1}{(2\alpha-1)^{2}(2\alpha^{2} - 1)} (2\alpha)^{2(n-p-1)} & \text{if $2\alpha^{2} > 1$}\end{cases}
\end{equation*}

It follows from \eqref{eq:cher-Scrit} that for all $\lambda > 0$, there exists a positive constant $c_{\alpha}$ such that
\begin{equation}\label{eq:cher-Scrit2}
\PP(b_{n}^{-1} R_{0}(n) > \delta) \leq \exp\left(-\lambda b_{n} |\GG_{n}|^{1/2} \, + \, c_{\alpha} \lambda^{2}  |\GG_{n-p}| \, + \, c_{\alpha} \lambda (1 + (2\alpha)^{n-p})\right).
\end{equation}
Taking $\lambda = 2^{-1} \, c_{\alpha}^{-1} \, \delta \, b_{n} \, |\GG_{p}| \, |\GG_{n}|^{-1/2}$ in \eqref{eq:cher-Scrit2}, we get, for some positive constant $c_{\alpha,\delta},$
\begin{equation*}\label{eq:cher-Scrit3}
\PP\left(b_{n}^{-1} R_{0}(n) > \delta\right) \leq \exp\left(-\frac{c_{\alpha,\delta} \, b_{n}^{2} \, |\GG_{p}| \, \delta^{2}}{4c_{\alpha}}\right).
\end{equation*}
Doing the same thing for the sequence $-\bF$ instead of $\bF$, we conclude that
\begin{equation*}\label{eq:cher-Scrit3}
\PP\left(\left|b_{n}^{-1} R_{0}(n)\right| > \delta\right) \leq 2 \, \exp\left(-\frac{c_{\alpha,\delta} \, b_{n}^{2} \, |\GG_{p}| \, \delta^{2}}{4c_{\alpha}}\right).
\end{equation*}
In the latter inequality, taking the $\log$, dividing by $b_{n}^{2}$ and letting $n$ goes to infinity, we get the result. 

\begin{rem}\label{rem:i-general}
Let $f \in \Bb_{b}(S)$. Since we will use frequently this type of inequality, we give here a general procedure to upper-bound the probability $\PP(||\GG_{n-p}|^{-1} M_{\GG_{n-p}}(\tilde{f})| > \delta).$ From Chernoff inequality, we have, for all $\lambda > 0$,
\begin{equation}\label{eq:i-C-gen}
\PP\left(\left||\GG_{n-p}|^{-1} M_{\GG_{n-p}}(\tilde{f})\right| > \delta\right) \leq \exp\left(- \lambda \delta |\GG_{n-p}|\right) \EE\left[\exp\left(\lambda M_{\GG_{n-p}}(\tilde{f})\right)\right].
\end{equation}
For all $m \in \{0, \ldots, n-p\}$, we set
\begin{equation*}
\II_{m} = \EE\left[\exp\left(2^{m} \lambda M_{\GG_{n-p-m}}(\Qq^{m} \tilde{f})\right)\right].
\end{equation*}
Using the branching Markov property, we have
\begin{equation*}
\II_{m} = \EE\left[\exp\left(2^{m+1} \lambda M_{\GG_{n-p-m-1}}(\Qq^{m+1} \tilde{f})\right) \JJ_{m}\right],
\end{equation*}
where
\begin{equation*}
\JJ_{m} = \prod_{i \in \GG_{n-p-m-1}} \EE_{X_{i}}\left[\exp\left(2^{m} \lambda \left(\Qq^{m}\tilde{f}(X_{i0}) + \Qq^{m}\tilde{f}(X_{i1}) - 2 \Qq^{m+1}\tilde{f}(X_{i})\right)\right)\right].
\end{equation*}
Using \eqref{eq:geom-erg} and Lemma \ref{lem:ABH}, we have the following upper-bound:
\begin{equation*}
\JJ_{m} \leq \exp\left(\lambda^{2} \|f\|^{2}_{\infty} M^{2} (1+\alpha)^{2} (2\alpha^{2})^{m} |\GG_{n-p}|\right).
\end{equation*}
This implies that
\begin{equation}\label{eq:i-gen-rec}
\II_{m} \leq \exp\left(\lambda^{2} \|f\|^{2}_{\infty} M^{2} (1+\alpha)^{2} (2\alpha^{2})^{m} |\GG_{n-p}|\right) \, \II_{m+1}.
\end{equation}
Using the recurrence relation \eqref{eq:i-gen-rec} and \eqref{eq:geom-erg} (to upper-bound $\II_{n-p}$), we are led to
\begin{equation}\label{eq:i-gen-rec2}
\II_{0} \leq \exp\left(\lambda^{2} \|f\|^{2}_{\infty} M^{2} (1+\alpha)^{2} a_{\alpha,n} |\GG_{n-p}| + \lambda \|f\|_{\infty} M (2\alpha)^{n-p}\right),
\end{equation}
where $a_{\alpha,n} = \sum_{m = 0}^{n-p-1} (2\alpha^{2})^{m}.$ We set $a_{\alpha} = \lim_{n \rightarrow \infty} a_{\alpha,n}$, which is finite since $2\alpha^{2} < 1.$ Taking $\lambda = \delta/(2\|f\|^{2}_{\infty} M^{2} (1 + \alpha)^{2} a_{\alpha})$ in \eqref{eq:i-C-gen} and using \eqref{eq:i-gen-rec2}, we are led to
\begin{equation*}
\PP\left(|\GG_{n-p}|^{-1} M_{\GG_{n-p}}(\tilde{f}) > \delta\right) \leq \exp\left(- \frac{\delta^{2}|\GG_{n-p}|}{4 \|f\|^{2}_{\infty}M^{2}(1 + \alpha)^{2} a_{\alpha}}\left(1 - \frac{2 \|f\|_{\infty} M \alpha^{n-p}}{\delta}\right)\right).
\end{equation*}
Finally, since we can do the same thing for $-f$ instead of $f$, we conclude that
\begin{equation}\label{eq:i-gen-end}
\PP\left(\left||\GG_{n-p}|^{-1} M_{\GG_{n-p}}(\tilde{f})\right| > \delta\right) \leq 2 \exp\left(- \frac{\delta^{2}|\GG_{n-p}|}{4 \|f\|^{2}_{\infty}M^{2}(1 + \alpha)^{2} a_{\alpha}}\left(1 - \frac{2 \|f\|_{\infty} M \alpha^{n-p}}{\delta}\right)\right).
\end{equation}
\end{rem}

\subsection{proof of \eqref{eq:cvsupexpR1}}\label{sec:mdpR1}
\quad\\
We set $g_{p} = \sum_{\ell = 0}^{p} 2^{p-\ell} \Qq^{p-\ell} \tilde{f}_{\ell}$ in such a way that using the definition of $R_{1}(n)$, we have
\begin{equation*}
\PP\left(b_{n}^{-1} \left|R_{1}(n)\right| > \delta\right) = \PP\left(|\GG_{n-p}|^{-1} \left|M_{\GG_{n-p}}\left(g_{p}\right)\right| > \delta b_{n} |\GG_{n}|^{-1/2} |\GG_{p}|\right).
\end{equation*}
Using \eqref{eq:geom-erg} and \eqref{eq:Bsupfell}, we have
\begin{equation*}
\|g_{p}\|_{\infty} \leq \begin{cases} c_{\infty} M (p+1) & \text{if $2\alpha \leq 1$} \\ c_{\alpha} c_{\infty} M & \text{if $1 < 2\alpha < \sqrt{2}$}.\end{cases}
\end{equation*}
Applying \eqref{eq:i-gen-end} to $g_{p}$ and $\delta b_{n} |\GG_{n}|^{-1/2} |\GG_{p}|$, we get, for $n$ going to infinity and for some positive constant $C_{\alpha,\delta}$,
\begin{equation*}
\PP\left(b_{n}^{-1} \left|R_{1}(n)\right| > \delta\right) \leq \begin{cases} 2 \exp\left(- C_{\alpha,\delta} \, \delta^{2} \, b_{n}^{2} \, |\GG_{p}|p^{-2}\right) & \text{if $2\alpha \leq 1$} \\ 2 \exp\left(- C_{\alpha,\delta} \, \delta^{2} b_{n}^{2} \, (2\alpha^{2})^{-p}\right) & \text{if $1 < 2\alpha < \sqrt{2}.$} \end{cases}
\end{equation*}
Finally, \eqref{eq:cvsupexpR1} follows by taking the $\log$, dividing by $b_{n}^{2}$ and letting $n$ goes to infinity in the latter inequality.

\subsection{Proof of \eqref{eq:mdp-Dnf-SCr}: Moderate deviations principle for $b_{n}^{-1}\Delta_n(\bF)$}\label{sec:mdpDelta}
\quad\\
First we study the bracket of $\Delta_n(\bF)$:
\begin{equation*}\label{eq:def-Vn-subc}
V(n)= \sum_{i\in \G_{n-p_n}} \E\left[ \Delta_{n, i}(\bF)^2|\cf_i\right]. 
\end{equation*}
Using \reff{eq:def-NiF} and \reff{eq:def-DiF}, we write:
\begin{equation}\label{eq:def-V}
V(n) = |\G_n|^{-1} \sum_{i\in \G_{n-p_n}} \E_{X_i}\left[\left(\sum_{\ell=0}^{p_n} M_{\G_{p_n-\ell}}(\tilde f_\ell) \right)^2 \right]-R_2(n)=V_1(n) +2V_2(n) - R_2( n),
\end{equation}
with:
\begin{align*}
V_1(n) & =   |\G_n|^{-1} \sum_{i\in \G_{n-p_n}} \sum_{\ell=0}^{p_n} \E_{X_i}\left[M_{\G_{p_n-\ell}}(\tilde f_\ell) ^2  \right] ,\\
V_2(n) & =  |\G_n|^{-1} \sum_{i\in \G_{n-p_n}} \sum_{0\leq \ell<k\leq p_n} \E_{X_i}\left[M_{\G_{p_n-\ell}}(\tilde f_\ell)  M_{\G_{p_n-k}}(\tilde f_k) \right], \\
R_2( n) &=\sum_{i\in \G_{n-p_n}} \E\left[ N_{n,i} (\bF) |X_i \right] ^2.
\end{align*}

\begin{lem}\label{lem:cvR2nScrit}
Under the Assumptions of Theorem \ref{theo:mdp-SC}, we have
\begin{equation}\label{eq:negR2nScrit}
\limsup_{n \rightarrow \infty} \frac{1}{b_{n}^{2}} \log\PP\left(R_{2}(n) > \delta\right) = -\infty.
\end{equation}
\end{lem}

\begin{proof}
Using the branching Markov property, we have
\begin{equation}\label{eq:def-R2n}
R_{2}(n) = |\GG_{n}|^{-1} \sum_{i \in \GG_{n-p}} g_{p}(X_{i}), \quad \text{with} \quad g_{p} = \left(\sum_{\ell = 0}^{p} 2^{p-\ell} \Qq^{p-\ell} \tilde{f}_{\ell}\right)^{2}.
\end{equation}
Using \eqref{eq:geom-erg} and \eqref{eq:Bsupfell}, we get
\begin{equation*}
\|g_{p}\|_{\infty} \leq c_{\infty}^{2} \, M^{2} \, \left(\sum_{\ell = 0}^{p} (2\alpha)^{p-\ell}\right)^{2} \leq \begin{cases} c_{\infty}^{2} \, M^{2} \, (p+1)^{2} & \text{if $2\alpha \leq 1$} \\ c_{\infty}^{2} \, M^{2} \, c_{\alpha} \, (2\alpha)^{2p}  & \text{if $1 < 2\alpha < \sqrt{2}$.}\end{cases} 
\end{equation*}
This implies that $R_{2}(n)$ is upper-bounded by a deterministic sequence which converge to $0$. As a consequence, we conclude that \eqref{eq:negR2nScrit} holds.
\end{proof}

\begin{lem}\label{cvV1nScrit}
Under the Assumptions of Theorem \ref{theo:mdp-SC}, we have
\begin{equation*}
\limsup_{n\rightarrow \infty } \frac{1}{b_{n}^{2}} \log \PP\left(\left|V_1(n) -\ssub_1(\bF)\right| > \delta\right) = -\infty,
\end{equation*}
where
\begin{equation*}\label{eq:S1}
\ssub_1(\bF)  =\sum_{\ell\geq 0} 2^{-\ell} \, \left\langle \mu,   \tilde f_\ell^ 2 \right\rangle + \sum_{\ell\geq 0, \, k\geq 0} 2^{k-\ell} \, \left\langle \mu, \cp\left((\cq^k \tilde f_\ell) \otimes^2\right)\right\rangle : = H_{3}(\bF) + H_{4}(\bF)
\end{equation*}
\end{lem}

\begin{proof}
Using \reff{eq:Q2}, we get:
\begin{equation}\label{eq:DV4nsub}
V_1(n)= V_3(n)+ V_4(n),
\end{equation}
with
\begin{align*}
V_3(n) &=  |\G_n|^{-1} \sum_{i\in \G_{n-p}} \sum_{\ell=0}^p 2^{p-\ell}\,\cq^{p-\ell} (\tilde f_\ell^2)(X_i),\\ 
V_4(n) &=     |\G_n|^{-1} \sum_{i\in \G_{n-p}} \sum_{\ell=0}^{p-1}\,\sum_{k=0}^{p-\ell -1} 2^{p-\ell+k} \,\cq^{p-1-(\ell+k)}\left(\cp\left(\cq^k \tilde f_\ell \otimes^2\right)\right)(X_i).  
\end{align*}
The proof is divided into two parts.

\medskip

\subsection*{Part I} First we prove that
\begin{equation}\label{eq:cvEV3nScrit}
\limsup_{n \rightarrow \infty} \frac{1}{b_{n}^{2}}\log\PP\left(|V_{3}(n) -  H_{3}(\bF)| > \delta\right) = -\infty.
\end{equation}
Since
\begin{equation*}
H_{3}(\bF) =  \sum_{\ell = 0}^{p} 2^{-\ell} \, \left\langle \mu,   \tilde f_\ell^ 2 \right\rangle +  \sum_{\ell > p} 2^{-\ell} \, \left\langle \mu,   \tilde f_\ell^ 2 \right\rangle \quad \text{and} \quad \lim_{n \rightarrow \infty} \sum_{\ell > p} 2^{-\ell} \, \left\langle \mu,   \tilde f_\ell^ 2 \right\rangle = 0 ,
\end{equation*}
then to get \eqref{eq:cvEV3nScrit}, it suffices to prove that
\begin{equation*}\label{eq:cvEV3nScrit2}
\limsup_{n \rightarrow \infty} \frac{1}{b_{n}^{2}}\log\PP\left(|V_{3}(n) -  H_{3}^{[n]}(\bF)| > \delta\right) = -\infty, \quad \text{where} \quad H_{3}^{[n]}(\bF) = \sum_{\ell = 0}^{p} 2^{-\ell} \, \left\langle \mu,   \tilde f_\ell^ 2 \right\rangle. 
\end{equation*}
We set
\begin{equation*}
g_{p} = \sum_{\ell = 0}^{p} 2^{-\ell} \Qq^{p-\ell}\left(\tilde{f}^{2}_{\ell} - \left\langle \mu,\tilde{f}^{2}_{\ell} \right\rangle\right) \quad \text{and then} \quad V_{3}(n)-H_{3}^{[n]}(\bF) = |\GG_{n-p}|^{-1} M_{\GG_{n-p}}(g_{p}).
\end{equation*}
Using \eqref{eq:geom-erg} and \eqref{eq:Bsupfell}, we have, for some positive constant $c_{\alpha}$,
\begin{equation*}
\|g_{p}\|_{\infty} \leq \begin{cases} 4 c_{\infty}^{2} c_{\alpha} M 2^{-p} & \text{if $2\alpha < 1$} \\ 4 c_{\infty}^{2} M (p+1) \alpha^{p} & \text{if $2\alpha >1$.}  \end{cases}  
\end{equation*}
Using \eqref{eq:i-gen-end}, we get, for $n$ going to infinity and for some positive constant $C_{\alpha,\delta}$:
\begin{equation*}
\PP\left(\left|V_{3}(n) -  H_{3}^{[n]}(\bF)\right| > \delta\right) \leq \begin{cases} 2 \exp\left(- \delta^{2} C_{\alpha,\delta} |\GG_{n+p}|\right) & \text{if $2\alpha \leq 1$} \\ 2\exp\left(-\delta^{2} C_{\alpha,\delta} p^{-2} |\GG_{n}| (2\alpha^{2})^{-p}\right) & \text{if $1 < 2\alpha < \sqrt{2}$.} \end{cases}
\end{equation*}
Finally, \eqref{eq:cvEV3nScrit} follows from the latter inequality by taking the $\log$ and dividing by $b_{n}^{2}$.
\medskip

\subsection*{Part II} Next, we prove that
\begin{equation}\label{eq:cvEV4nScrit}
\limsup_{n \rightarrow \infty} \frac{1}{b_{n}^{2}}\log\PP\left(|V_{4}(n) -  H_{4}(\bF)| > 2\delta\right) = -\infty.
\end{equation}
Note that $V_{4}(n) - H_{4}(\bF) = |\GG_{n-p}|^{-1} M_{\GG_{n-p}}(H_{4,n}(\bF) - H_{4}(\bF))$, where
\begin{equation}\label{eq:def-H4nf}
H_{4,n}(\bF) = \sum_{\ell \geq 0, k \geq 0} h_{\ell,k}^{(n)} \ind_{\{\ell + k < p\}} \quad \text{and} \quad H_{4}(\bF) = \sum_{\ell \geq 0, k \geq 0} h_{\ell,k}, 
\end{equation}
with
\begin{equation*}
h_{\ell,k}^{(n)} = 2^{k-\ell}\,  \cq^{p-1-(\ell+k)} \left(\cp\left(\cq^k \tilde f_\ell \otimes^2 \right)\right) \quad \text{and} \quad h_{\ell,k} = \left\langle \mu, \cp\left((\cq^k \tilde f_\ell) \otimes^2\right)\right\rangle.
\end{equation*}
Using \eqref{eq:geom-erg} and \eqref{eq:Bsupfell}, we have
\begin{equation}\label{eq:ih+hlkn-Sc}
|h_{\ell,k}| + |h_{\ell,k}^{(n)}| \leq 2 C_{\infty}^{2} M^{2} (2\alpha^{2})^{k} 2^{-\ell}.
\end{equation}
Let $r_{0}$ large enough such that
\begin{equation}\label{eq:i-delta-Sc}
2 c_{\infty}^{2} M^{2} \sum_{\ell \vee k > r_{0}} (2\alpha^{2})^{k} 2^{-\ell} \leq \delta.
\end{equation}
For $n$ going to infinity, we have
\begin{align}
|M_{\GG_{n-p}}(H_{4,n}(\bF) - H_{4}(\bF))| &\leq |M_{\GG_{n-p}}(\sum_{\ell \vee k \leq r_{0}} (h_{\ell,k}^{(n)} - h_{\ell,k}))| + M_{\GG_{n-p}}(\sum_{\ell \vee k > r_{0}} (|h_{\ell,k}| + |h_{\ell,k}^{(n)}|)) \nonumber \\
& \leq |M_{\GG_{n-p}}(\sum_{\ell \vee k \leq r_{0}} (h_{\ell,k}^{(n)} - h_{\ell,k}))| \, + \, 2 \, c_{\infty}^{2} \, M^{2} |\GG_{n-p}| \, \sum_{\ell \vee k > r_{0}} (2\alpha^{2})^{k} 2^{-\ell}, \label{eq:i-H4n-H4fSc} 
\end{align}
where we used \eqref{eq:ih+hlkn-Sc} for the second inequality. From \eqref{eq:i-H4n-H4fSc}, we get 
\begin{equation}\label{eq:i-H4-V4Sc}
|V_{4}(n) - H_{4}(\bF)| \leq  |\GG_{n-p}|^{-1} |M_{\GG_{n-p}}(g_{p})| \, + \, 2 \, c_{\infty}^{2} \, M^{2}\, \sum_{\ell \vee k > r_{0}} (2\alpha^{2})^{k} 2^{-\ell},
\end{equation}
where $g_{p} = \sum_{\ell \vee k \leq r_{0}} (h_{\ell,k}^{(n)} - h_{\ell,k}).$
From \eqref{eq:i-delta-Sc} and \eqref{eq:i-H4-V4Sc}, to get \eqref{eq:cvEV4nScrit}, it suffices to prove that
\begin{equation}\label{eq:cvEV4nScrit2}
\limsup_{n \rightarrow \infty} \frac{1}{b_{n}^{2}} \log \PP\left(\left| |\GG_{n-p}|^{-1} M_{\GG_{n-p}}(g_{p})\right| > \delta \right).
\end{equation}
Using \eqref{eq:geom-erg} and \eqref{eq:Bsupfell} twice, we have, for some positive constant $c_{\alpha}$,
\begin{equation*}
\|g_{p}\|_{\infty} \leq c_{\infty}^{3} M^{3} c_{\alpha} \gamma(r_{0}) \alpha^{p-1} \quad \text{where} \quad \gamma(r_{0}) = \begin{cases} r_{0} (2\alpha)^{r_{0}} & \text{if $2\alpha \leq 1$} \\ (2\alpha)^{r_{0}} & \text{if $1 < 2\alpha < \sqrt{2}$}.\end{cases}
\end{equation*}
Using \eqref{eq:i-gen-end} with $g_{p}$ instead of $f$, we get, for some positive constant $C_{\alpha,\delta}$,
\begin{equation*}
\PP\left(\left| |\GG_{n-p}|^{-1} M_{\GG_{n-p}}(g_{p})\right| > \delta \right) \leq \exp\left(- \delta^{2} C_{\alpha,\delta} |\GG_{n-p}| \alpha^{-2p}\right).
\end{equation*}
Taking the $\log$, dividing by $b_{n}^{2}$ and letting $n$ goes to infinity in the latter inequality, we get \eqref{eq:cvEV4nScrit2} and then \eqref{eq:cvEV4nScrit}.
\end{proof}

\begin{lem}\label{cvV2nScrit}
Under the Assumptions of Theorem \ref{theo:mdp-SC}, we have
\begin{equation*}
\limsup_{n\rightarrow \infty } \frac{1}{b_{n}^{2}} \log \PP\left(\left|V_2(n) -\ssub_2(\bF)\right| > \delta\right) = -\infty,
\end{equation*}
where
\begin{equation*}\label{eq:S2}
\ssub_2(\bF) = \sum_{0\leq \ell< k} 2^{-\ell} \left\langle \mu,   \tilde f_k \cq^{k-\ell} \tilde f_\ell\right\rangle + \sum_{\substack{0\leq \ell< k\\ r\geq 0}} 2^{r-\ell} \left\langle \mu, \cp\left( \cq^r \tilde f_k \sot \cq^{k-\ell+r} \tilde f_\ell  \right)\right\rangle : = H_{5}(\bF) + H_{6}(\bF)
\end{equation*}
\end{lem}

\begin{proof}
Using \reff{eq:Q2-bis}, we get:
\begin{equation}\label{eq:decom-V2}
V_2(n) = V_5(n)+ V_6(n),
\end{equation}
with
\begin{align*}
V_5(n) &=  |\G_n|^{-1} \sum_{i\in \G_{n-p}} \sum_{0\leq \ell<k\leq  p } 2^{p-\ell} \cq^{p-k} \left( \tilde f_k \cq^{k-\ell} \tilde f_\ell\right)(X_i),\\
V_6(n) &=     |\G_n|^{-1} \sum_{i\in \G_{n-p}} \sum_{0\leq \ell<k<  p }\sum_{r=0}^{p-k-1}  2^{p-\ell+r} \, \cq^{p-1-(r+k)}\left(\cp\left(\cq^r \tilde f_k \sot \cq  ^{k-\ell+r} \tilde f_\ell \right)\right)(X_i).
\end{align*}
\subsection*{Part I} First, we prove that
\begin{equation}\label{eq:cvEV6nScrit}
\limsup_{n \rightarrow \infty} \frac{1}{b_{n}^{2}}\log\PP\left(|V_{6}(n) -  H_{6}(\bF)| > 2\delta\right) = -\infty.
\end{equation}
We have $V_{6}(n) - H_{6}(\bF) = |\GG_{n-p}|^{-1} M_{\GG_{n-p}}(H_{6,n}(\bF) - H_{6}(\bF))$, where
\begin{equation*}\label{eq:def-Hn-V6}
H_{6,n}(\bF) = \sum_{\substack{0\leq \ell< k \\ r\geq 0}}  h_{k,\ell,r}^{(n)}\,\ind_{\{r+k<  p\}} \quad \text{and} \quad H_{6}(\bF) = \sum_{\substack{0\leq \ell< k \\ r\geq 0}} h_{k, \ell,r},
\end{equation*}
with
\begin{multline*}
h_{k, \ell,r}^{(n)} =  2^{r-\ell} \,\cq^{p-1-(r+k)}\left(\cp\left(\cq^r \tilde f_k \sot \cq  ^{k-\ell+r} \tilde f_\ell \right)\right) \quad \text{and} \\ 
h_{k,\ell,r} = 2^{r-\ell} \left\langle \mu, \cp\left(\cq^r \tilde f_k \sot \cq^{k-\ell+r} \tilde f_\ell\right)\right\rangle.
\end{multline*}
Using \eqref{eq:geom-erg} and \eqref{eq:Bsupfell}, we have
\begin{equation}\label{eq:ih+hlkrn-Sc}
|h_{\ell,k,r}| + |h_{\ell,k,r}^{(n)}| \leq 2 \, c_{\infty}^{2} \, M^{2} \, 2^{r-\ell} \, \alpha^{k - \ell + 2r}.
\end{equation}
Let $r_{0}$ large enough such that
\begin{equation}\label{eq:eq:ideltaSV6}
2 \, c_{\infty}^{2} \, M^{2} \, \sum_{\substack{0 \leq \ell < k \\ r \geq 0 \\ k \vee r > r_{0}}}  2^{r-\ell} \, \alpha^{k - \ell + 2r} < \delta.
\end{equation}
We set $g_{p} = \sum_{0\leq \ell< k, \, \, r \geq 0, \, \, k \vee r \leq r_{0}} (h^{(n)}_{k,\ell,r} - h_{k,\ell,r}).$ Using \eqref{eq:ih+hlkrn-Sc}, we have, for $n$ going to infinity in such a way that $p > r_{0}$, 
\begin{equation}\label{eq:i-H6-V6Sc}
\left|V_{6}(n) - H_{6}(\bF)\right| \leq |\GG_{n-p}| \left|M_{\GG_{n-p}}(g_{p})\right| + 2 \, c_{\infty}^{2} \, M^{2} \, \sum_{\substack{0 \leq \ell < k \\ r \geq 0 \\ k \vee r > r_{0}}}  2^{r-\ell} \, \alpha^{k - \ell + 2r}.
\end{equation}
From \eqref{eq:eq:ideltaSV6} and \eqref{eq:i-H6-V6Sc}, to get \eqref{eq:cvEV6nScrit}, it suffices to prove that
\begin{equation}\label{eq:cvEV6nScrit2}
\limsup_{n \rightarrow \infty} \frac{1}{b_{n}^{2}} \log \PP\left(\left| \, |\GG_{n-p}|^{-1} M_{\GG_{n-p}}\left(g_{p}\right)\right| > \delta \right) = -\infty.
\end{equation}
Using \eqref{eq:geom-erg} and \eqref{eq:Bsupfell} twice, we have, for some positive constant $c_{\alpha}$,
\begin{equation*}
\|g_{p}\|_{\infty} \leq 2 \, c_{\alpha} \, c_{\infty}^{3} \, M^{3} \, \gamma(r_{0}) \, \alpha^{p},
\end{equation*}
where
\begin{equation*}
\gamma(r_{0}) = \begin{cases} (2\alpha)^{-r_{0}} & \text{if $2\alpha < 1$} \\ r_{0}^{2} & \text{if $2\alpha \geq 1$.} \end{cases}
\end{equation*}
Using \eqref{eq:i-gen-end} with $g_{p}$ instead of $f$, we get, for some positive constant $C_{\alpha,\delta}$,
\begin{equation*}
\PP\left(\left| |\GG_{n-p}|^{-1} M_{\GG_{n-p}}(g_{p})\right| > \delta \right) \leq \exp\left(- \delta^{2} C_{\alpha,\delta} |\GG_{n-p}| \alpha^{-2p}\right).
\end{equation*}
Taking the $\log$, dividing by $b_{n}^{2}$ and letting $n$ goes to infinity in the latter inequality, we get \eqref{eq:cvEV6nScrit2} and then \eqref{eq:cvEV6nScrit}.

\subsection*{Part II} Next, with the finite constant $H_5(\bF)$ defined by:
\begin{equation*}\label{eq:def-H5}
H_5(\bF) = \sum_{0\leq \ell< k}
2^{-\ell} \langle \mu,   \tilde f_k \cq^{k-\ell} \tilde f_\ell\rangle,
\end{equation*} 
we prove that
\begin{equation}\label{eq:cvEV5nScrit}
\limsup_{n \rightarrow \infty} \frac{1}{b_{n}^{2}}\log\PP\left(|V_{5}(n) -  H_{5}(\bF)| > 2\delta\right) = -\infty.
\end{equation}
We set
\begin{equation*}
H_{5,n}(\bF) = \sum_{0\leq \ell < k} h_{k,\ell}^{(n)} \ind_{\{k \leq p\}}, \quad \text{and} \quad H_{5}^{[n]}(\bF) = \sum_{0\leq \ell < k} h_{k,\ell} \ind_{\{k \leq p\}},
\end{equation*}
with
\begin{equation*}
h_{k,\ell}^{(n)} = 2^{-\ell} \Qq^{p-k}\left(\tilde{f}_{k} \Qq^{k-\ell} \tilde{f}_{\ell}\right) \ind_{\left\{k \leq p\right\}} \quad \text{and} \quad h_{\ell,k} = \left\langle \mu, \tilde{f}_{k} \Qq^{k-\ell} \tilde{f}_{\ell}\right\rangle.
\end{equation*}
We have the following decomposition: 
\begin{equation*}
V_{5}(n) -  H_{5}(\bF) = |\GG_{n-p}|^{-1} M_{\GG_{n-p}}\left(H_{5,n}(\bF) -  H_{5}^{[n]}(\bF)\right) + \left(H_{5}^{[n]}(\bF) - H_{5}(\bF)\right). 
\end{equation*}
Using \eqref{eq:geom-erg} and \eqref{eq:Bsupfell}, we have
\begin{equation*}
|h_{k,\ell}^{(n)}| + |h_{k,\ell}| \leq 2 c_{\infty}^{2} M \alpha^{k-\ell} 2^{-\ell}, 
\end{equation*} 
which implies that $\lim_{n \rightarrow \infty} |H_{5}(\bF) - H_{5}^{[n]}(\bF)| = 0.$ As a result, to get \eqref{eq:cvEV5nScrit}, it suffices to prove that
\begin{equation}\label{eq:cvEV5nScrit2}
\limsup_{n \rightarrow \infty} \frac{1}{b_{n}^{2}}\log\PP\left(\left||\GG_{n-p}|^{-1} M_{\GG_{n-p}}\left(H_{5,n}(\bF) -  H_{5}^{[n]}(\bF)\right)\right| > \delta\right) = -\infty.
\end{equation}
Setting $g_{p} = H_{5,n}(\bF) -  H_{5}^{[n]}(\bF)$, we have, using \eqref{eq:geom-erg} and \eqref{eq:Bsupfell}:
\begin{equation*}
\|g_{p}\|_{\infty} \leq \begin{cases} c_{\alpha} 2^{-p} & \text{if $2\alpha < 1$} \\ c_{\alpha} p \alpha^{p} & \text{if $1 \leq 2\alpha < \sqrt{2}$,}\end{cases}
\end{equation*}
for some positive constant $c_{\alpha}.$ Finally, \eqref{eq:cvEV5nScrit2}, and then \eqref{eq:cvEV5nScrit}, follows by applying \eqref{eq:i-gen-end} to $g_{p}$ instead of $f$ and by taking the $\log$, dividing by $b_{n}^{2}$ and by letting $n$ goes to infinity.
\end{proof}
As a direct consequence of \eqref{eq:def-V} and Lemmas \ref{lem:cvR2nScrit}, \ref{cvV1nScrit} and \ref{cvV2nScrit}, we have the following result.
\begin{lem}\label{lem:cvbrac-SC}
Under the Assumptions of Theorem \ref{theo:mdp-SC}, we have
 \begin{equation*}\label{eq:cvBra-C}
\limsup_{n \rightarrow \infty} \frac{1}{b_n^2}\log \PP\left( \vert V(n)- \ssub (\bF)\vert > \delta \right) = -\infty
\end{equation*}
\end{lem}

We can now state the following result.
\begin{lem}
Under Assumptions of Theorem \ref{theo:mdp-SC}, we have that $b_{n}^{-1} \Delta_{n}(\bF)$ satisfies a moderate deviation principle with speed $b_{n}^{2}$ and rate function $I$ defined in \eqref{eq:mdp-rate}.
\end{lem}
\begin{proof}
Since $p < n/2$, we have for all $i \in \GG_{n-p}$,
\begin{equation*}
|\Delta_{n,i}(\bF)| \leq 2 c_{\infty} 2^{-\tfrac{n}{2} + p} \leq C,
\end{equation*}
where $C$ is a positive constant. This implies that $\Delta_{n}(\bF)$ is a martingale with bounded differences. Using the result of Dembo \cite{dembo1996moderate} (see also Djellout \cite{djellout2002moderate} and Puhalski \cite{puhalskii1997large}), we get the result from Lemma \ref{lem:cvbrac-SC}.
\end{proof}

\subsection{Completion of the proof of Theorem \ref{theo:mdp-SC}.}
Finally, using the decomposition \eqref{eq:N=D+R} and the results of sections \ref{sec:supER0}, \ref{sec:mdpR1} and \ref{sec:mdpDelta}, we deduce Theorem \ref{theo:mdp-SC}.

\medskip

\section{Proof of Theorem \ref{theo:mdp-C}}\label{sec:proof-T2}

\subsection{A quick overview of our strategy}
\quad\\
Let $(p_n, n\in \N)$ be a non-decreasing sequence of elements of $\N^*$ such that, for all $\lambda>0$:
\begin{equation}\label{eq:def-pn}
p_n< n, \quad \lim_{n\rightarrow \infty } p_n/n=1 \quad \text{and}\quad \lim_{n\rightarrow \infty } n-p_n - \lambda \log(n)=+\infty .
\end{equation}
When there is no ambiguity, we write $p$ for $p_n$. Let us  consider the sequence $\bF=\left(f_{\ell}, \ell \in \mathbb{N} \right)$ of elements of  $\mathcal{B}_b(S)$ which satisfies the Assumption  \eqref{hyp:F3} (and then Assumption \ref{hyp:F2}) uniformly, namely:
\begin{equation}\label{hyp:F4}
\vert \mathcal{Q}^n(\tilde{f}_\ell) \vert \leq M \alpha^n c_{\infty} \quad \text{and} \quad \vert \mathcal{Q}^n(\hat{f}_\ell) \vert \leq M \beta_n \alpha^{n} c_{\infty} 
\end{equation}
It follows from \eqref{hyp:F4} that  there exists a finite constant    $c_J$ depending only on $ \{\alpha_j, j \in J \}$ such that for all $\ell \in \NN $, $n \in \NN$, $j_0 \in J $
\begin{equation*}\label{hyp:F5}
\vert f_{\ell} \vert  \leq Mc_{\infty},\quad \vert \tilde{f}_{\ell} \vert  \leq Mc_{\infty},\quad \vert  \left \langle\mu, f_{\ell} \right\rangle \vert  \leq Mc_{\infty},\quad |\sum_{j \in J}\theta_j^n\crr_{j}(f_{\ell}) | \leq 2Mc_{\infty} \quad \text{and}\quad  \vert \crr_{{j}_0}(f_{\ell}) \vert \leq  c_JMc_{\infty}.
\end{equation*}
We recall that:
\begin{equation}\label{eq:N=D+R}
N_{n, \emptyset}(\bF) = \Delta_n(\bF) +R_0(n)+R_1(n),
\end{equation}
with
\begin{equation*}\label{eq:reste01}
R_0(n)= |\G_n|^{-1/2}\sum_{k=0}^{n-p_n-1} M_{\G_k}(\tilde f_{n-k}) \quad\text{and}\quad R_1(n)= \sum_{i\in \G_{n-p_n}}\E\left[N_{n,i}(\bF) |\, \cf_i\right].
\end{equation*}
Let $(b_n)_{n \in \mathbb{N}}$ be a sequence elements of $\mathbb{N}$ such that :
\begin{equation*}\label{eq:vitesse-S}
b_n \rightarrow \infty \qquad  \text{and } \qquad \frac{b_n}{\sqrt{n \vert \mathbb{G}_{n}\vert}}\underset{n   \to \infty}{\longrightarrow} 0
\end{equation*}
Our goals will be achieved if we prove the following:
\begin{align}
&\label{eq:cvsupexpR01}
\forall \delta > 0, \quad \lim_{n \rightarrow \infty} \frac{1}{b_{n}^{2}} \log \PP(|b_{n}^{-1}n^{-1/2} R_{0}(n)| > \delta) = - \infty; \\
&\label{eq:cvsupexpR11}
\forall \delta > 0, \quad \lim_{n \rightarrow \infty} \frac{1}{b_{n}^{2}} \log \PP(|b_{n}^{-1}n^{-1/2} R_{1}(n)| > \delta) = - \infty; \\
&\label{eq:mdp-Dnf-C}
b_{n}^{-1}n^{-1/2}\Delta_n(\bF) \quad \text{satisfies a MDP on $S$ with speed $b_n^2$ and rate function $I$.}
\end{align}
\subsection{proof of \eqref{eq:cvsupexpR01}}
\quad\\
We follow the same lines of  the proof of \eqref{eq:cvsupexpR0} with $2\alpha^2=1$. First, using Chernoff inequality, we have 
\begin{equation*}\label{eq:devr0}
\mathbb{P}\left( b_n^{-1}n^{-1/2}R_0(n) > \delta\right) \leq \exp\left(- \lambda b_n\delta \vert \mathbb{G}_n \vert^{\frac{1}{2}}n^{1/2}\right) \mathbb{E}\left[\exp\left(\lambda \sum_{k=0}^{n-p-1} M_{\mathbb{G}_k}(\tilde{f}_{n-k})\right)\right].
\end{equation*}
Next, taking $\lambda =\tfrac{b_n \delta  \left(n \vert \mathbb{G}_{n} \vert\right)^{1/2} }{2 c_{\alpha}(n-p)\vert \mathbb{G}_{n-p} \vert}$ and doing the same thing for $-\bF$ instead of $\bF$, we get
\begin{equation*}
\mathbb{P}\left(\vert b_n^{-1}n^{-1/2} R_0(n)  \vert >\delta\right)   \vert \leq 2\exp\left( -\frac{b_n^2 \delta^2 n \vert \mathbb{G}_p \vert }{4 c_{\alpha}(n-p)} \right).
\end{equation*}
Finally, taking the $\log$, dividing by $b_n^2$ and letting $n$ goes to infinity, we get the result.

 \medskip
 
\begin{rem}\label{rem:C-gen}
We have the following version of Remark \ref{rem:i-general} when $2\alpha^{2} = 1$:
 \begin{equation}\label{rem2}
\mathbb{P}\left( \vert \vert \mathbb{G}_{n-p} \vert^{-1} M_{\mathbb{G}_{n-p}}(\tilde{f}) \vert > \delta n^{1/2}    \right) \leq 2 \exp\left( -\frac{\delta^2 n \vert \mathbb{G}_{n-p}\vert}{4 \Vert f \Vert_{\infty}^2M^2(1+\alpha)^2(n-p)} \right).
 \end{equation}
 \end{rem}
 
 \subsection{proof of \eqref{eq:cvsupexpR11}}
With $g_p= \sum_{\ell = 0}^{p} 2^{p - \ell}\mathcal{Q}^{p - \ell}\tilde{f}_{\ell}$, and  using the definition of $R_{1}(n)$, we have for all $\delta > 0$
\begin{align*}
\mathbb{P}\left(\vert  b_n^{-1}R_1(n) \vert  >\delta n^{1/2}\right) =\mathbb{P}\left( \vert \mathbb{G}_{n-p} \vert ^{-1} \vert  M_{\mathbb{G}_{n-p}(g_p)} \vert  >b_n\delta n^{1/2} \vert \mathbb{G}_p \vert\vert \mathbb{G}_n \vert^{-1/2} \right). 
\end{align*}

So according to \eqref{hyp:F4} , we have:
\begin{equation*}\label{gg}
\Vert g_p \Vert_{\infty} \leq c_{\infty}c_{\alpha} M \vert \mathbb{G}_{p} \vert^{1/2}.
\end{equation*}

By applying $\eqref{rem2}$ to $g_p$ and $ b_n\delta  n^{1/2}\vert \mathbb{G}_p \vert\vert \mathbb{G}_n \vert^{-1/2} $  and using the fact that  $2\alpha^2=1$, we have:
\begin{align*}
\mathbb{P}\left(\vert  b_n^{-1}R_1(n) \vert  >\delta n^{1/2}  \right) \leq 2 \exp\left(-\frac{b_n^2 \delta^2 n}{4c_{\infty}^2c_{\alpha}^2M^4(1+\alpha)^2(n-p)}\right).
\end{align*}
So taking the $\log$ and dividing by $b_n^2$, and letting n goes to infinity, we get the result.

\medskip

\subsection{Proof of \eqref{eq:mdp-Dnf-C}: Moderate deviations principle for  $b_n^{-1}n^{-1/2}\Delta _n(\bF)$}
\quad\\
First we study the bracket of $n^{-\frac{1}{2}}\Delta_n(\bF)$ given by $n^{-1} V(n),$ where $V(n)$ is defined in \eqref{eq:def-V}.  We have the following result:

\begin{lem}\label{cv-expR2-C}
Under the assumptions of Theorem \ref{theo:mdp-C}, we have
\begin{equation*}\label{23}
\underset{n \to \infty}{\limsup}\frac{1}{b_n^2}\log \PP\left( \vert n^{-1}R_2(n) \vert > \delta \right) = - \infty.
\end{equation*} 
\end{lem}
\begin{proof}
Recall the definition of $R_{2}(n)$ and $g_{p}$ given in \eqref{eq:def-R2n}. So according to \eqref{hyp:F4} and using $2\alpha^2 = 1$, we have $\Vert g_p \Vert_{\infty} \leq c_{\infty}^2c_{\alpha}^2 M^2\vert \mathbb{G}_{p} \vert.$ This implies that
\begin{equation*}
R_2(n) \leq \frac{c_{\infty}^2c_{\alpha}^2 M^2 }{n} \underset{n \to \infty}{\longrightarrow} 0.
\end{equation*}
Therefore, $R_2(n)$ is upper-bounded by a deterministic sequence which goes to $0$. According to Remark \ref{rem3}, we get the result.
\end{proof}

\begin{lem}\label{cv-expoV1-crit}
Under the assumptions of Theorem \ref{theo:mdp-C}, we have
\begin{equation*}\label{21}
\limsup_{n \rightarrow \infty}\frac{1}{b_n^2}\log \PP\left( \vert n^{-1}V_1(n)- \scrit_1 (\bF)\vert > \delta \right)=-\infty.
\end{equation*}
\end{lem}
\begin{proof}
Recall the decomposition of $V_{1}(n)$ given in \eqref{eq:DV4nsub} and  the definition of $\scrit_{1}(\bF)$ given in \ref{eq:S1-crit}. The proof is divided into two parts.
\subsection*{Part I}
First we prove that
\begin{equation*}
\limsup_{n \rightarrow \infty} \frac{1}{b_{n}^{2}}\log\PP\left(|n^{-1}V_{3}(n)| > \delta \right) = -\infty.
\end{equation*}
Indeed we have 
\begin{equation*}
n^{-1}V_3(n) = \vert \G_{n-p} \vert^{-1}n^{-1} M_{\G_{n-p}}(g_p), \quad \text{where} \quad g_{p} = \sum_{\ell = 0}^{p} 2^{-\ell}\mathcal{Q}^{p - \ell}(\tilde{f}_{\ell}^{2}).
\end{equation*}
Since the sequence $\bF = (f_{\ell}, \ell \in \NN)$ satisfies \eqref{eq:Bsupfell},  we have $\Vert g_p \Vert_{\infty} \leq 4 c_{\infty}^{2}.$ This implies that $|n^{-1}V_3(n)| \leq 4 c_{\infty}^{2} n^{-1}.$ So $n^{-1}V_3(n)$ is upper-bounded by a deterministic sequence which goes to 0. Then applying the remark \ref{rem3} to $n^{-1}V_3(n)$, we get the result.
\subsection*{Part II}
Next, we prove that
\begin{equation*}
\limsup_{n \rightarrow \infty} \frac{1}{b_{n}^{2}}\log\PP\left(|n^{-1}V_{4}(n) -  \scrit_{1}(\bF)| > \delta \right) = -\infty.
\end{equation*}
Recall $H_{4,n}(\bF)$ and $f_{k,\ell}$ defined respectively in \eqref{eq:def-H4nf} and \eqref{eq:def-fkl*}. For $k, \ell,r\in \N$, we consider the following $\C\text{-}$valued functions defined on $S^{2}$: 
\begin{equation}\label{eq:def-fkl}
f_{k,\ell,r} = \Big( \sum_{j\in J} \theta_j^{r}\,  \crr_j (f_k) \Big) \sot \Big(\sum_{j\in J}  \theta_j^{r+k-\ell}\,  \crr_j (f_\ell) \Big).
\end{equation}
Recall that $2\alpha^{2} = 1.$ We set $\bar H_{4,n} = \sum_{\ell\geq 0, \, k\geq 0} \bar{h}_{\ell,k}^{(n)}\,  \ind_{\{\ell+k<  p\}}$ with
\begin{equation*}
\bar h_{\ell,  k}^{(n)} = 2^{k-\ell}\,\alpha^{2k}\,   \cq^{p-1-(\ell+k)} \left(\cp f_{\ell,\ell,k}\right) = 2^{-\ell}\,   \cq^{p-1-(\ell+k)} \left(\cp f_{\ell,\ell,k}\right).
\end{equation*}
For $f \in \Bb_{b}(S)$, recall $\hat{f}$ defined in \eqref{eq:tilde-hat-f}. Then we have $h_{\ell,k}^{(n)} - \bar{h}_{k,\ell}^{(n)} = h_{\ell,k}^{n,1} + h_{\ell,k}^{n,2} + h_{\ell,k}^{n,3}$, where
\begin{align*}
&h_{\ell,k}^{n,1} = 2^{k-\ell} \Qq^{p-1-(\ell + k)} \Pp(\Qq^{k} \hat{f}_{k} \sot \Qq^{k} \hat{f}_{k}),\\
&h_{\ell,k}^{n,2} = 2^{k-\ell} \Qq^{p-1-(\ell + k)}\Pp(\Qq^{k} \hat{f}_{k} \sot \Qq^{k} (\sum_{j \in J} \crr_{j}(f_{k}))),\\
&h_{\ell,k}^{n,3} = 2^{k-\ell} \Qq^{p-1-(\ell + k)}\Pp(\Qq^{k} (\sum_{j \in J} \crr_{j}(f_{k})) \sot \Qq^{k} (\sum_{j \in J} \crr_{j}(f_{k}))).
\end{align*}
This implies that
\begin{equation*}
n^{-1} |\GG_{n-p}|^{-1} M_{\GG_{n-p}}(H_{4,n}(\bF) - \bar{H}_{4,n}(\bF)) = n^{-1} |\GG_{n-p}|^{-1} \sum_{u \in \GG_{n-p}} \sum_{i=1}^{3} \sum_{\ell \geq 0; k \geq 0} h_{\ell,k}^{n,i}(X_{u}) \ind_{\{\ell + k < p\}}.
\end{equation*}
Using Assumption \ref{hyp:F3}, \eqref{eq:Bsupfell} and the fact that the sequence $(\beta_{k}, k \in \NN)$ is decreasing, we can upper bound each function $|h_{\ell,k}^{n,i}|$, $i \in \{1,2,3\}$, by $C 2^{-\ell} \beta_{k}.$ This implies that
\begin{equation}\label{eq:B-H4n-H4bn}
\left|n^{-1} |\GG_{n-p}|^{-1} M_{\GG_{n-p}}(H_{4,n}(\bF) - \bar{H}_{4,n}(\bF))\right| \leq C n^{-1} \sum_{k = 0}^{p-1} \beta_{k}.
\end{equation}
 We set 
 \begin{equation*}\label{eq:H4starn}
 H_{4}^{[n]} = \sum_{\ell \geq 0;k\geq 0} h_{\ell,k} \ind_{\{k+\ell<p\}} \quad \text{with} \quad h_{\ell,k} = 2^{-\ell} \langle \mu,\Pp(f_{\ell,\ell,k}) \rangle.
 \end{equation*}
 Using Assumption \ref{hyp:F3} and \eqref{eq:Bsupfell}, we get
 \begin{equation*}
 \left|\bar h_{\ell,  k}^{(n)} - h_{\ell,k}\right| \leq 2^{-\ell} \left(\frac{1}{\sqrt{2}}\right)^{p-1-(\ell+k)} \left\|\Pp(f_{\ell,\ell,k})\right\|_{\infty} \leq C \left(\frac{1}{\sqrt{2}}\right)^{p-1} \left(\frac{1}{\sqrt{2}}\right)^{\ell-k}. 
 \end{equation*}
 This implies that
 \begin{equation}\label{eq:B-Hbarn4-H4n}
 \left|n^{-1} |\GG_{n-p}|^{-1} M_{\GG_{n-p}}(\bar{H}_{4,n}(\bF) - H_{4}^{[n]}(\bF))\right| \leq C n^{-1}.
 \end{equation}
 Finally, from \cite{BD2020}, we have
 \begin{equation}\label{eq:L-H4n-sC1}
 \lim_{n \rightarrow \infty} |n^{-1}H_{4}^{[n]}(\bF) - \scrit_{1}(\bF)| = 0.
 \end{equation}
From \eqref{eq:B-H4n-H4bn}, \eqref{eq:B-Hbarn4-H4n} and \eqref{eq:L-H4n-sC1}, we conclude that  $\vert n^{-1}V_{4}(n) - \scrit_{1}(\bF) \vert $ is upper-bounded by a deterministic sequence which goes to $ 0$. Therefore applying the remark \ref{rem3} to $n^{-1}V_{4}(n) - \scrit_{1}(\bF) $, we get the result.
\end{proof}

\begin{lem}\label{cv-expoV2-crit}
Under the assumptions of Theorem \ref{theo:mdp-C}, we have
\begin{equation*}\label{22}
\limsup_{n \rightarrow \infty}\frac{1}{b_n^2}\log \PP\left( \vert n^{-1}V_2(n)- \scrit_2 (\bF)\vert > \delta \right)=-\infty.
\end{equation*}
\end{lem}
\begin{proof}
Recall the decomposition given in \eqref{eq:decom-V2}. Then following the lines of the proof of Lemma \ref{cv-expoV1-crit}, we prove that
\begin{equation*}
\limsup_{n \rightarrow \infty} \frac{1}{b_{n}^{2}}\log\PP\left(|n^{-1}V_{5}(n)| > \delta \right) = -\infty \, \,  \text{and} \, \, \limsup_{n \rightarrow \infty} \frac{1}{b_{n}^{2}}\log\PP\left(|n^{-1}V_{6}(n) -  \scrit_{2}(\bF)| > \delta \right) = -\infty,
\end{equation*}
and the result follows.
\end{proof}
As a direct consequence of \eqref{eq:def-V} and Lemmas \ref{cv-expR2-C}, \ref{cv-expoV1-crit} and \ref{cv-expoV2-crit}, we have the following result.
\begin{lem}\label{lem:cvbracket-C}
Under the assumptions of Theorem \ref{theo:mdp-C}, we have
 \begin{equation*}\label{eq:cvBra-C}
\limsup_{n \rightarrow \infty} \frac{1}{b_n^2}\log \PP\left( \vert V(n)- \scrit (\bF)\vert > \delta \right) = -\infty.
\end{equation*}
\end{lem}

Next, contrary to the sub-critical case, we need to check the exponential Lindeberg condition and Chen-Ledoux type condition, that is conditions \textbf{(C2)} and \textbf{(C3)} given in Proposition \ref{prop:mdp}. Indeed, in the critical case, the martingale $n^{-\frac{1}{2}}\Delta_{n}(\bF)$ does not have bounded differences in such a way that Lemma \ref{lem:cvbracket-C} is not longer sufficient to get the moderate deviations principle of $n^{-\frac{1}{2}} \Delta_{n}(\bF)$. In order to check exponential Lindeberg condition, we have the following exponential Lyapunov condition which implies the  exponential Lindeberg condition.
\begin{lem}\label{lem:lyapunov-C}
Under the assumptions of Theorem \ref{theo:mdp-C}, we have
\begin{equation*}\label{eq:lyapunov}
\limsup_{n \rightarrow \infty} \frac{1}{b_{n}^{2}} \log \PP\left( \sum_{i \in \GG_{n-p}} n^{-2} \EE\left[\Delta_{n,i}(\bF)^{4}|\Ff_{i}\right] > \frac{\delta |\GG_{n}| n}{b_{n}^{2}} \right) = -\infty \quad \forall \delta > 0.
\end{equation*}
\end{lem}
\begin{proof}
For all $i \in \GG_{n-p}$, we have
\begin{equation}\label{eq:ED4iSMG}
\EE\left[\Delta_{n,i}(\bF)^{4}|\Ff_{i}\right] \leq 16 (p+1)^{3} 2^{-2n} \sum_{\ell = 0}^{p} \EE_{X_{i}}\left[M_{\GG_{p-\ell}}(\tilde{f}_{\ell})^{4}\right],
\end{equation}
where we have  used the definition of $\Delta_{n,i}(\bF)$, the inequality $(\sum_{k=0}^r a_k)^4 \leq  (r+1)^3 \sum_{k=0}^r a_k^4$ and the branching Markov property. Using \eqref{eq:geom-erg} and \eqref{eq:Bsupfell}, we can apply Theorem 2.1 given in \cite{BDG14} with $2 \alpha^{2} = 1$ to get $\EE_{X_{i}}[M_{\GG_{p-\ell}}(\tilde{f}_{\ell})^{4}] \leq C p^{2} 2^{2(p-\ell)}.$ The latter inequality and \eqref{eq:ED4iSMG} imply that
\begin{equation}\label{eq:ineg-lya}
\frac{b_{n}^{2}}{n |\GG_{n}|}\sum_{i \in \GG_{n-p}} n^{-2} \EE\left[\Delta_{n,i}(\bF)^{4}|\Ff_{i}\right] \leq C n^{3} 2^{-n + p} ((n|\GG_{n}|)^{-1} b_{n}^{2}).
\end{equation}
From \eqref{eq:vitesse-Sc} and \eqref{eq:def-pn}, we have
\begin{equation}\label{eq:lim-n3npnGb}
\lim_{n \rightarrow \infty} n^{3} 2^{-n + p} ((n|\GG_{n}|)^{-1} b_{n}^{2}) = 0.
\end{equation}
Finally, the result of the Lemma follows using \eqref{eq:ineg-lya}, \eqref{eq:lim-n3npnGb} and Remark \ref{rem3}.  
\end{proof}

For Chen-Ledoux type condition, we have the following result.
\begin{lem}\label{lem:chen-ledoux}
Under the assumptions of Theorem \ref{theo:mdp-C}, we have
\begin{equation*}\label{eq:chen-L}
\limsup_{n \rightarrow \infty} \frac{1}{b_{n}^{2}} \log \PP\left(|\GG_{n}| \sup_{i \in \GG_{n-p}} \PP_{\Ff_{i}}\left(|\Delta_{n,i}(\bF)| > b_{n} \sqrt{n |\GG_{n}|}\right)\right) = - \infty.
\end{equation*}
\end{lem}

\begin{proof}
	\quad\\
Using \eqref{eq:Bsupfell}, we have, for $n$ large enough and for all $i \in \GG_{n-p}$, $|\Delta_{n,i}(\bF)| \leq C 2^{-\tfrac{n}{2} + p} < b_{n}\sqrt{n|\GG_{n}|}$.\\ This implies that
\begin{equation*}
\PP_{\Ff_{i}}\left(|\Delta_{n,i}(\bF)| > b_{n} \sqrt{n |\GG_{n}|}\right) = 0 \quad \forall i \in \GG_{n-p}.
\end{equation*}
From the latter equality and using the convention $\log(0) = - \infty,$ we get the result of the Lemma.
\end{proof}
We can now state the following result.
\begin{lem}\label{lem:mdp-D-Sc}
Under the assumptions of Theorem \ref{theo:mdp-C}, we have that $n^{-1/2} b^{-1}_{n} \Delta_{n}(\bF)$ satisfies a moderate deviation principle with speed $b_{n}^{2}$ and rate function $I$ defined in \eqref{eq:mdp-rate}. 
\end{lem}
\begin{proof}
Applying Theorem 1 in \cite{djellout2002moderate} (a simplified version is given in Proposition \ref{prop:mdp}) to the martingale differences $n^{-1/2} \Delta_{n,i}(\bF)$, the proof follows  from Lemmas \ref{lem:cvbracket-C}, \ref{lem:lyapunov-C} and \ref{lem:chen-ledoux}.  
\end{proof}

\subsection{Completion of the proof of Theorem \ref{theo:mdp-C}.}
Finally, using \eqref{eq:N=D+R}, \eqref{eq:cvsupexpR01}, \eqref{eq:cvsupexpR11} and Lemma \ref{lem:mdp-D-Sc}, we deduce Theorem \ref{theo:mdp-C}.

\section{Numerical studies}\label{sec:numerical}  

\begin{figure}
	\centering
	\begin{subfigure}{0.45\textwidth} 
		\includegraphics[width=\textwidth]{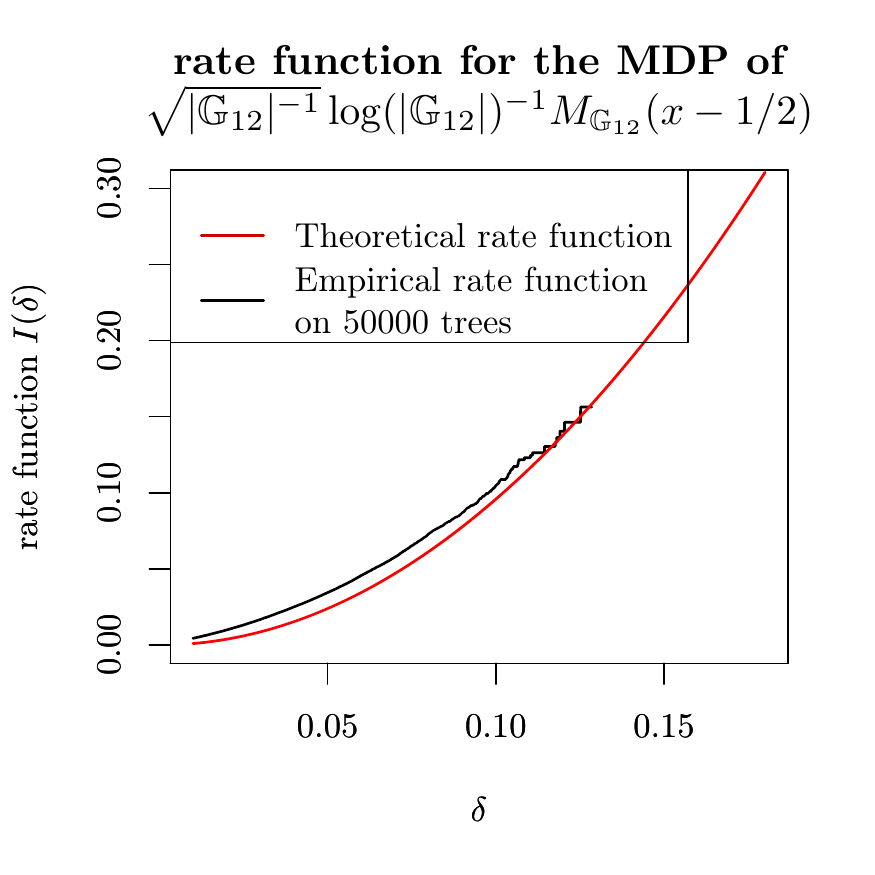}
	\end{subfigure}
	\begin{subfigure}{0.45\textwidth} 
		\includegraphics[width=\textwidth]{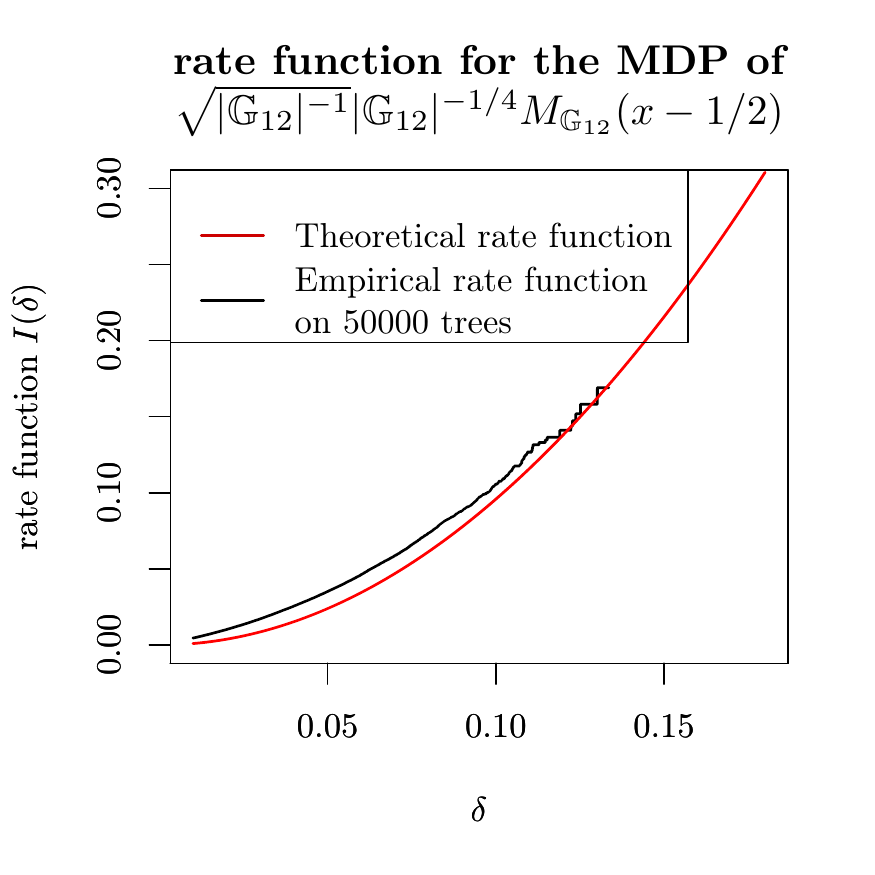}
	\end{subfigure}
	\centering
	\begin{subfigure}{0.45\textwidth} 
		\includegraphics[width=\textwidth]{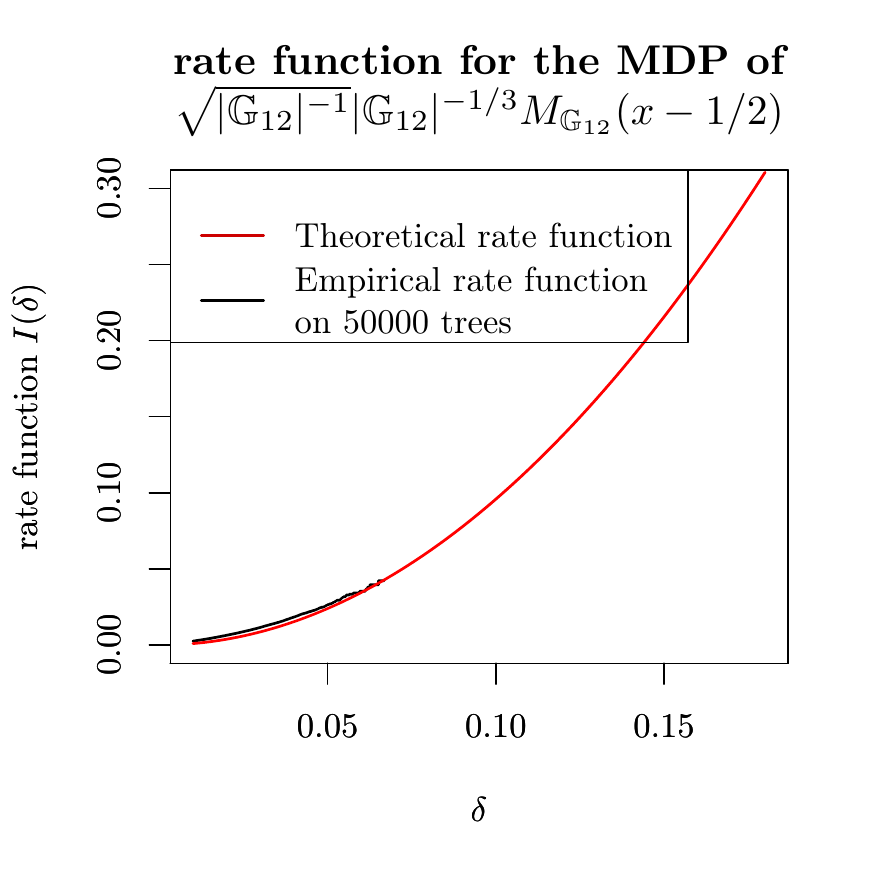}
	\end{subfigure}
	\begin{subfigure}{0.45\textwidth} 
		\includegraphics[width=\textwidth]{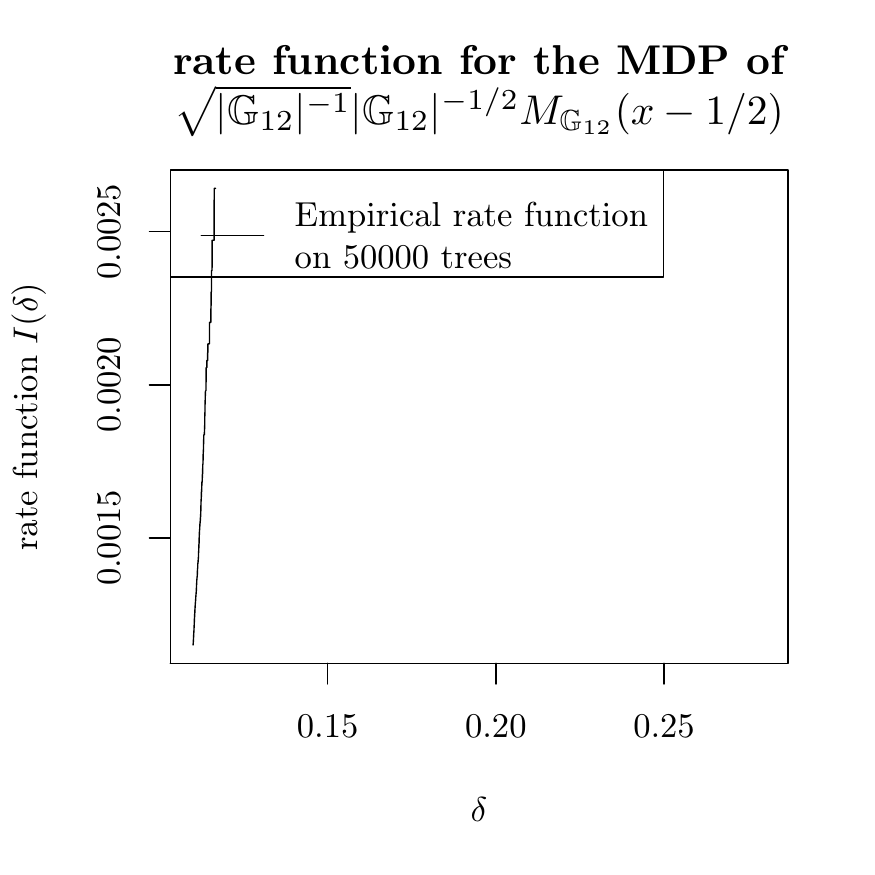}
	\end{subfigure}
        \caption{Exact and empirical rate functions for the moderate deviation principles of $|\GG_{n}|^{-1/2}b_{n}^{-1}M_{\GG_{n}}(x-1/2).$ In the first three figures, one can see that empirical rate function fit well the exact rate function. The differences can be explained from the fact that the sample size is not large enough to generate enough large deviation events. In the last figure, one can see that the empirical rate function is reduced to $0$. This is due to the fact that the speed considered here is valid only in the critical case, not in the subcritical case.} \label{fig:mdp-rate}
\end{figure}

For our numerical illustrations, we consider a BMC $(X_{u}, u \in \TT)$ living in $[0,1]$, with transition $\Pp = \Qq \otimes \Qq$ given by 
\begin{equation*} \label{BAR transition}
\Qq(x,y) : = (1 - x) \, \frac{y(1-y)^{2}}{B(2,3)} + x \, \frac{y^{2}(1-y)}{B(3,2)} \, , \quad x, y  \in [0,1],
\end{equation*}
with $B(\alpha, \beta)$ the normalizing constant of a standard Beta distribution with shape parameters $\alpha$ and $\beta$. For simplicity, we choose $X_{\emptyset}$ such that $\Ll(X_{\emptyset}) = Beta(2,2)$, where $Beta(2,2)$ is the standard Beta distribution with shape parameters $(2,2)$. Now, one can prove that this process is stationary, it has an explicit invariant density: the standard Beta distribution with shape parameters $(2,2)$. One can also prove that $ \EE\left[X_{u0}|X_{u}\right] = \EE\left[X_{u1}|X_{u}\right] =  X_{u}/5 + 2/5,$  (for more details, we refer e.g. to \cite{pitt2002constructing}). Now, it is not hard to verify that this process satisfies our required assumptions. In particular, using for example Theorem 2.1 in \cite{hairer2011yet}, one can prove that Assumption \ref{hyp:F2} is satisfied with $\alpha = 1/5.$ We are thus in the sub-critical case. First, we will illustrate Theorem \ref{theo:mdp-SC} (and more precisely Corollary \ref{cor:mdp-SC}) with the sequence $\bF = (f, 0, 0, \ldots)$ and the function $f(x) = x.$ In this case, we have the following exact results:
\begin{equation*}
\langle \mu,f \rangle = \frac{1}{2}, \quad \Qq^{k}\tilde{f}(x) = 5^{-k}( x - \tfrac{1}{2} ) \quad \forall k \geq 0, \quad \Sigma_{\GG}(f) = 6/115 \quad \text{and} \quad I(\delta) = \frac{115}{12}\delta^{2}.
\end{equation*}
Next, we will illustrate that the range of speed considered in the critical case does not work in this example. For that purpose, we simulate $B = 50000$ samples $(X^{(s)} = (X^{(s)}_{u},u \in \GG_{12}), s \in \{1,\ldots, B\})$ of the bifurcating Markov chain at the $n\text{-}$th generation, with $n = 12$. For each sample $X^{(s)}$, we compute $b_{n}^{-1} N_{n,\emptyset}^{(s)}(f) = b_{n}^{-1}|\GG_{n}|^{-1/2} \sum_{u \in \GG_{n}} (X_{u}^{(s)} - 1/2).$ Finally, for different values of $\delta > 0$, we compute $b_{n}^{-2} \log(B^{-1} \sum_{s = 1}^{B} \ind_{\{|b_{n}^{-1} N_{n,\emptyset}^{(s)}(f)| > \delta\}}).$  This allows us to get empirical values of the rate function. Next in the same graph, we plot the true rate function and the empirical rate function. As we can see in Figure \ref{fig:mdp-rate}, the empirical rate function fit well exact rate function, except in the last figure where the empirical rate function is near to $0$ since the speed considered is not valid for the subcritical case, but only for the critical. We also stress that the differences observed between empirical and exact rate functions can be explained from the fact that the sample size is not large enough.

\section{Appendix}\label{sec:appendix}
	
The following Remark is used in  the proofs of Lemmas \ref{21}, \ref{22} and \ref{23}.
\begin{rem}\label{rem3}
\color{black}
\quad\\
We assume that $(S,d)$ is a metric space.
Let  $\left(Z_n\right)_{n \in \NN}$ be a sequence of random variables valued in $S$,  $Z$ a random variable valued in $S$ and $v_n$ a rate. So if  $d(Z_n,Z) $ is  upper-bounded by  a deterministic  sequence which converges to  $0$,  then, for all sequence $(v_{n}, n \in \NN)$ converging to $\infty$, $Z_n$  converges $ v_{n}$- superexponentially  fast in probability to $Z$, that is for all $\delta > 0$,
\begin{equation*}
\limsup_{n \rightarrow \infty} \frac{1}{v_n} \log\PP \left( d(Z_n,Z) >\delta\right) =-\infty.
\end{equation*}
\end{rem}
\medskip

The following result is known as Azuma-Bennett-Hoeffding inequality
\cite{Azuma, Bennett, Hoeffding}.

\begin{lem}\label{lem:ABH}
Let $X$ be a real-valued and centered random variable such that
$\displaystyle a\leq X\leq b$ a.s., with $a<b$. Then for all
$\lambda>0$, we have
\begin{equation*}
\E\left[\exp\left(\lambda X\right)\right] \leq
\exp\left(\frac{\lambda^{2}(b-a)^{2}}{8}\right).
\end{equation*}
\end{lem}
We have the following many-to-one formulas. Ideas of the proofs can be found in \cite{Guyon} and \cite{BDG14}.
\begin{lem}\label{lem:Qi}
Let $f,g\in \cb(S)$, $x\in S$ and $n\geq m\geq 0$. Assuming that all the quantities below are well defined,  we have:
\begin{align}
\label{eq:Q1} \E_x\left[M_{\G_n}(f)\right] &=|\G_n|\, \cq^n f(x)= 2^n\, \cq^n f(x) ,\\
\label{eq:Q2} \E_x\left[M_{\G_n}(f)^2\right] &=2^n\, \cq^n (f^2) (x) + \sum_{k=0}^{n-1} 2^{n+k}\,   \cq^{n-k-1}\left( \cp \left(\cq^{k}f\otimes \cq^k f \right)\right) (x),\\
\label{eq:Q2-bis} \E_x\left[M_{\G_n}(f)M_{\G_m}(g)\right] &=2^{n} \cq^{m} \left(g \cq^{n-m} f\right)(x)\\
\nonumber &\hspace{2cm} + \sum_{k=0}^{m-1} 2^{n+k}\, \cq^{m-k-1} \left(\cp\left(\cq^k g \sot \cq^{n-m+k} f\right) \right)(x). 
\end{align}
\end{lem}

We recall here a simplified version of Theorem 1 in \cite{djellout2002moderate}. We consider the real martingale $(M_{n},n\in\NN)$ with respect to the filtration $(\Hh_{n}, n\in \NN)$ and we denote $(\langle M\rangle_{n}, n \in \NN)$ its bracket. 

\begin{prop}\label{prop:mdp} Let $(b_{n})$ a sequence satisfying
\[b_{n} \quad \text{is increasing}, \quad b_{n} \longrightarrow +\infty, \quad \frac{b_{n}}{\sqrt{n}}\longrightarrow 0,\]
such that $c(n) := \sqrt{n}/b_{n}$ is non-decreasing, and define the
reciprocal function $c^{-1}(t)$ by
\[c^{-1}(t):=\inf\{n\in \mathbb{N}: c(n)\geq t\}.\]
Under the following conditions:
\begin{enumerate}
\item [(\textbf{C1})] there exists $Q\in \mathbb{R}_{+}^{*}$ such that for all $\delta > 0$,

$\displaystyle \limsup_{n \rightarrow \infty} \frac{1}{b_{n}^{2}} \log\left(\PP\left(\left|\frac{\langle M\rangle_{n}}{n} - Q\right| > \delta\right)\right) = -\infty,$

\item [(\textbf{C2})]  $\displaystyle \limsup\limits_{n\rightarrow
+\infty}\frac{1}{b_{n}^{2}}\log\left(n \quad \underset{1\leq k\leq
c^{-1}(b_{n+1})}{\rm ess\,sup}
\mathbb{P}(|M_{k}-M_{k-1}| > b_{n}\sqrt{n} \Big|\mathcal{H}_{k-1})\right)=-\infty,$

\item [(\textbf{C3})] for all $ a>0$ and for all $\delta > 0$,

$\displaystyle \limsup_{n \rightarrow \infty} \frac{1}{b_{n}^{2}} \log\left(\PP\left(\frac{1}{n}\sum\limits_{k=1}^{n}\mathbb{E}\left(|M_{k}-M_{k-1}|^{2} \mathbf{1}_{\{|M_{k}-M_{k-1}| \geq a \frac{n}{b_{n}}\}} \Big|\mathcal{H}_{k-1} \right) > \delta\right)\right) = - \infty,$
\end{enumerate}
$\left(M_{n}/(b_{n}\sqrt{n})\right)_{n\in \N}$ satisfies the MDP in $\mathbb{R}$ with the
speed $b_{n}^{2}/n$ and the rate function $\displaystyle I(x) =\frac{x^{2}}{2Q}.$
\end{prop}

\begin{rem}
For all $n \geq 1$, we set $m_{n} = M_{n} - M_{n-1}$. Note that, in Proposition \ref{prop:mdp}, if the sequence $(m_{n})_{n \geq 1}$ is uniformly bounded, we recover a simplified version of the result of Dembo \cite{dembo1996moderate} and if the sequence $(m_{n})_{n \geq 1}$ is bounded by a deterministic sequence, we recover a simplified version of the result of Puhalskii \cite{puhalskii1997large}. 
\end{rem}

\bibliographystyle{abbrv}
\bibliography{biblio}
\end{document}